\documentclass[10pt]{IEEEtran}
\usepackage{graphicx}
\usepackage[cmex10]{amsmath}
\usepackage{amsthm}

\usepackage{mdwmath}
\usepackage{mdwtab}

\usepackage{color}

\usepackage{eqparbox}

\usepackage{url}

\newtheorem{thm}{Theorem}
\newtheorem{proposition}[thm]{Proposition}
\newtheorem{lemma}[thm]{Lemma}
\theoremstyle{remark}
\newtheorem{remark}{Remark}
\theoremstyle{definition}
\newtheorem{definition}{Definition}
\newtheorem{assumption}{Assumption}

\begin{document}
\title{Weakly-coupled Systems in Quantum Control}

\author{Nabile~Boussa\"{i}d,
        Marco~Caponigro,
        and~Thomas~Chambrion
\thanks{N. Boussa\"{i}d is with the Universit\'e de Franche--Comt\'e,
Laboratoire de math\'ematiques, 16 route de Gray, 25030
Besan\c{c}on Cedex, France.}
\thanks{M. Caponigro is with \'Equipe M2N, Conservatoire National des Arts et M\'etiers, 292, rue Saint-Martin, 75003, Paris, France.}%
\thanks{T. Chambrion is with the Universit\'e de Lorraine, Institut {\'E}lie Cartan de Nancy, 172 boulevard des Aiguillettes, BP 70239, 54506 Vand{\oe}uvre Cedex, France and with
INRIA, 54600 Villers-l\`es-Nancy, France.}
\thanks{This  work has been supported by the INRIA Nancy-Grand Est ``CUPIDSE'' Color
program.
Second and third authors were partially supported by French Agence National de
la Recherche ANR ``GCM'', program ``BLANC-CSD'', contract number NT09-504590.
The third author was partially supported by European Research Council ERC StG
2009 ``GeCoMethods'', contract number 239748. 
}
}
\maketitle

\begin{abstract}
Weakly-coupled systems are a class of infinite dimensional conservative
bilinear control systems with discrete spectrum. An important feature of these systems is that they can be
precisely approached by finite dimensional Galerkin approximations. This
property is of particular interest for the approximation of quantum system
dynamics and the control of the bilinear Schr\"odinger equation.

The present study provides rigorous definitions
and analysis of the dynamics of
weakly-coupled systems and gives sufficient conditions for an infinite
dimensional quantum control system to be weakly-coupled.  As an illustration we
provide examples chosen among common physical systems.
\end{abstract}



%

\section{Introduction}

\subsection{Physical context}
The state of a quantum system evolving on a finite
dimensional
Riemannian manifold $\Omega$, with associated measure $\mu$, is described by
its \emph{wave function}, that is, an element of the unit sphere of $L^2(\Omega,
\mathbf{C})$. Any physical quantity $\mathcal O$ (e.g. energy, position, momentum) is associated with a Hermitian operator $O:L^2(\Omega,
\mathbf{C}) \rightarrow L^2(\Omega,
\mathbf{C})$. The expected value of $\mathcal O$ for a system with wave function $\psi$ is  equal to $\displaystyle{\int_{\Omega} \!\overline{\psi(x)} O\psi(x)
\mathrm{d}\mu(x)}$.   For instance, a system with wave function $\psi$ is in a subset $\omega$ of
$\Omega$ with probability $\displaystyle{\int_{\omega} \!|\psi|^2
\mathrm{d}\mu}$.

The dynamics of a closed system submitted to excitations by $p$ external  fields ({\it e.g.} lasers) is described, under the dipolar approximation, by the bilinear
Schr\"odinger
equation 
\begin{eqnarray}
\mathrm{i}\frac{\partial}{\partial t} \psi(x,t)&=&-\frac{1}{2}\Delta \psi +V(x)
\psi(x,t) \nonumber\\
& & \quad
+ \sum_{l=1}^p u_l(t)
W_l(x) \psi(x,t),\label{eq:blse}
\end{eqnarray}
where $\Delta$ is the Laplace-Beltrami operator on $\Omega$,
$V:\Omega\rightarrow \mathbf{R}$ is a real function, usually called potential,
carrying the physical properties of the uncontrolled system,
$W_l:\Omega\rightarrow \mathbf{R}$, $1\leq l \leq p$, is a real function
modeling a laser $l$, and $u_l$, $1\leq l \leq p$, usually
called control, is a real function of the time representing the intensity of
the laser $l$.

In recent years there has been an increasing interest in studying the
controllability of the bilinear Schr\"odinger equation~\eqref{eq:blse} mainly 
due  to its importance for many advanced applications such as Nuclear Magnetic
Resonance, laser spectroscopy, and quantum information science. The problem
concerns the existence of control laws $(u_1,\ldots,u_p)$ steering the system
from a given initial state to a pre-assigned final state in a given time. 
Considerable efforts have been made  to study this problem and the main
difficulty is the fact that the state space, namely $L^2(\Omega,\mathbf{C})$, has infinite dimension. Indeed in~\cite{bms}, a result which implies (see
\cite{turinici}) strong limitations on the exact controllability of
the bilinear Schr\"odinger Equation has been proved. Hence, one has to look for
weaker controllability properties as, for instance, approximate controllability
or controllability between eigenstates of the Sch\"odinger operator (which are
the most relevant cases from the physical viewpoint).
In dimension one, in the case $p=1$, and for a specific class of control
potentials a description of the reachable set has been provided~\cite{beauchard,
camillo}.
In dimension larger than one or in  more general situations, the exact
description of the reachable set appears to be more difficult and at the moment
only approximate controllability results are available (see for
example~\cite{Nersy, Schrod, Schrod2} and references therein).

\subsection{Finite dimensional approximations}

 To avoid difficulties in dealing with infinite dimensional systems,
for instance in practical computations or simulations, one can project
system (\ref{eq:blse}) on finite dimensional subspaces of
$L^2(\Omega,\mathbf{C})$. 
A vast literature is currently available on control of
bilinear finite dimensional quantum systems (see for
instance~\cite{dalessandro-book} and references therein)
 thanks, also, to general controllability methods for left-invariant control systems on compact Lie group~\cite{book}.
 A crucial issue is to guarantee that the finite
dimensional approximations have dynamics close to the one of
the original infinite dimensional system.

{In \cite{Schrod} and \cite{Schrod2}, precise
estimates
of the distance
between the infinite dimensional systems and some of its Galerkin approximations
are used to prove that
systems of type (\ref{eq:blse}) are approximately controllable under physical conditions of 
non-degeneracy of the discrete spectrum of $-\Delta +V$. These estimates are derived for
a sequence of \emph{ad hoc} controls designed to steer the system from a given source to a
given target. Besides the discretness of the spectrum of $-\Delta+V$, very few regularity 
assumptions are made on (\ref{eq:blse}). Since the potential $W$ is not assumed to be
bounded or regular (say, not even continuous), the estimates obtained for
a control $u$ can possibly fail to hold for 
controls close to
$u$, for instance, in a small neighborhood of $u$ for some $H^k$ norm. 

Such pathological irregularities (everywhere discontinuous potentials or wave
functions) are physically irrelevant. Following~\cite[Chapter 2.A]{cohen77},
real potentials are at least continuous and
wave functions are smooth (i.e., infinitely differentiable). As a consequence,
most of the potentials and wave functions encountered in the literature 
are analytic. This strong regularity allows stronger estimates than those in \cite{Schrod} and \cite{Schrod2}.

As a matter of fact, a special class of bilinear systems of the type of
(\ref{eq:blse}), called weakly-coupled (see Definition~\ref{DEF_weakly_coupled} in Section
\ref{SEC_weakly_coupled}), exhibits very nice properties of approximations (see
Theorem~\ref{PRO_Good_Galerkin_approximation} below). {Physically,
for such weakly-coupled systems, the energy $\int_{\Omega}\overline{(-\Delta +V)\psi} ~\psi
\mathrm{d}\mu$ is bounded by an explicit function of the $L^1$ norm of the control $u$ (see Proposition~\ref{PRO_croissance_norme_A}), preventing propagation of  the wave function to  high energy levels.   

The notion of weakly-coupled systems, and the fact that such systems can be precisely approached by finite dimensional bilinear systems, has many applications.

First, taking advantage of the powerful tools of the geometric control theory
for finite dimensional systems \cite{dalessandro-book},
 this definition can be used for the analysis and the open-loop control of infinite dimensional bilinear quantum systems. 
For instance, we used this method  to prove that the rotational wave approximation (which is classical for 
finite dimensional systems) is still valid for infinite dimensional systems (\cite{ACCFEPS}) or to exhibit 
an example of bilinear system approximately controllable in arbitrary small times (\cite{Time}).

Second, it provides easily computable bounds on the size of the finite dimensional systems to be considered for the numerical simulations of systems of type (\ref{eq:blse}) in order to guarantee a given upper bound for the error.
This has been used in \cite{ACCgates} to implement a quantum gate in an infinite
dimensional systems modelling the rotation of a 2D-molecule. 

While the notion of weakly-coupled systems has been originally developed for
open-loop control, the approximation results apply without modification (both for the theoretical
analysis and the numerical simulation) in broader contexts. An example of Lyapunov design of open-loop   control is presented in Section  \ref{SEC_Lyapunoiv}.
}

The aim of this work is to provide an analysis of weakly-coupled systems, to present a sufficient condition for controllability for these systems, and to
show that two important types of bilinear quantum systems frequently encountered
in the literature are weakly-coupled.

\subsection{Content of the paper}
In Section~\ref{SEC_weakly_coupled} we introduce the notion of weakly-coupled systems for bilinear quantum systems and we state some properties of their finite dimensional approximations. In particular the most important property of this class of systems is that they have a Good Galerkin Approximation (Theorem~\ref{PRO_Good_Galerkin_approximation}), that is, a finite dimensional approximation whose dynamics is arbitrarily close the the one of the infinite dimensional system. Thanks to this feature we are able to show an approximate controllability result in higher norm for such a class of systems (Proposition~\ref{prop:controlHs}).

We then
study two important examples of weakly-coupled systems, the first (Section
\ref{SEC_compact_case}) covering, among others, the case where $\Omega$ is
compact (Section~\ref{sec:compact})
and the second (Section~\ref{SEC_tridiagonal}) the case where the system \eqref{eq:blse}
is tri-diagonal.

\section{Weakly-coupled systems}\label{SEC_weakly_coupled}

\subsection{Abstract framework}
We reformulate~\eqref{eq:blse} in a more abstract framework using
the language of functional analysis.  This
reformulation allows us  to treat examples slightly more general
than \eqref{eq:blse}, for instance, the example in  Section~\ref{SEC_trap_ion}. 
 For the convenience of the reader, we recall some 
basic notions of operator theory in the appendix.

In a separable Hilbert
space $H$
endowed with norm $\| \cdot \|$ and Hilbert product $\langle \cdot, \cdot
\rangle$, we consider the evolution problem
\begin{equation}\label{eq:main}
\frac{\mathrm{d} \psi}{\mathrm{d}t}= \left(A+ \sum_{l=1}^pu_l B_l \right)\psi
\end{equation}
where $(A,B_1,\ldots,B_p)$ satisfies Assumption~\ref{ass:ass}.

\begin{assumption}\label{ass:ass}
$(A,B_1,\ldots,B_p)$ is a $(p+1)$-uple of linear operators such that
 \begin{enumerate}
  \item for every $u$ in $\mathbf{R}^{p}$, $A+\sum_{l}u_{l}B_{l}$ is essentially
skew-adjoint on the domain $D(A)$ of $A$ and $\mathrm{i}(A+\sum_{l}u_{l}B_{l})$ is bounded from
below;
\label{EgaliteDomaine}
  \item $A$ is skew-adjoint and has purely discrete spectrum $(-\mathrm{i}
\lambda_j)_{j \in \mathbf{N}}$, the sequence  $(\lambda_j)_{j \in \mathbf{N}}$
is {positive non-decreasing and unbounded.}
\label{ASS_bounded_from_below}
 \end{enumerate}
\end{assumption}

In the rest of our study, we denote by $(\phi_j)_{j \in
\mathbf{N}}$ a Hilbert basis of $H$ such
that $A\phi_j=-\mathrm{i}  \lambda_j \phi_j$ for every $j$ in $\mathbf{N}$.
We denote by $D(A+\sum_{l}u_{l}B_{l})$ the domain where
$A+\sum_{l}u_{l}B_{l}$ is skew-adjoint.

Assumption~\ref{ass:ass}.\ref{EgaliteDomaine} ensures that, for every
constants $u_1,\ldots,u_p $ in $\mathbf{R}$, $A+\sum_l u_l B_l$ generates a
group of unitary propagators. Hence, for every initial state
$\psi_0$ in
$H$, for every piecewise constant control $u:t\in\mathbf{R}\rightarrow
\sum_{n=0}^N u^n\chi_{[t_n,t_{n+1})}(t)\in\mathbf{R}^p$,  where
$\chi_{[a,b)}(t)$ stands for the characteristic function of the interval $[a,b)$, with $0=t_0\leq t_1
\leq
\ldots  \leq t_{N+1}$ we can define the
solution of \eqref{eq:main} by $t\mapsto \Upsilon^u_t
\psi_0$, where
\begin{eqnarray*}
\lefteqn{\Upsilon^u_t=e^{(t-t_{j-1})(A +\sum u_l^{j-1} B_l)}\circ}\\
&& e^{(t_{j-1}-t_{j-2})(A
+\sum u_l^{j-2} B_l)} \circ \cdots  \circ e^{t_{0}(A+\sum u_l^0  B_l)},
\end{eqnarray*}
for $t\in [t_{j-1},t_{j})$, $j=1,\ldots,N$.

\begin{remark}
 From Assumption~\ref{ass:ass}.\ref{EgaliteDomaine} we deduce
that the resolvent of $A$ is compact, and for every
$u\in{\mathbf R}^p$, $A+\sum_l u_lB_l$ is bounded from $D(A)$ to $H$ as well as
$\sum_l u_lB_l$. As a consequence,
the Resolvent Identity~\eqref{eq:ResolventIdentity} applied to $A+\sum_l u_lB_l$ and $A$, gives that
the
resolvent of $A+\sum_l u_lB_l$ is compact.
The spectrum of $-\mathrm{i}(A+\sum_l u_lB_l)$ is discrete and accumulates only
at $+ \infty$ (as $-\mathrm{i}(A+\sum_l u_lB_l)$ is bounded	
from below).
\end{remark}

\subsection{Energy growth}
From Assumption~\ref{ass:ass}.\ref{ASS_bounded_from_below}, the operator
$A$ is self-adjoint with positive eigenvalues. For every $\psi$ in
$D(A)$, $\mathrm{i}A\psi=\sum_{j\in \mathbf{N}} \lambda_j \langle \phi_j, \psi
\rangle \phi_j$. For every $s\geq 0$,
using~\eqref{eq:popopo}
we define the $s$-norm by $\|\psi \|_s=\||A|^s \psi\|$
for every $\psi$ in $D(|A|^s)$. The $1/2$-norm plays an important role in
physics: for every $\psi$ in $D(|A|^{1/2})$, the quantity $|\langle A \psi,\psi
\rangle |=\|\psi\|_{1/2}^2$ is the expected value of the energy.
In the case $s=0$, we have the Hilbert space norm, thus we write $\|\psi\|$
instead of $\|\psi\|_0$.

\begin{remark}
The $s$-norm is a way to measure the regularity of the wave
functions. In the case where $|A|$ is the Laplace-Beltrami operator of
a smooth compact manifold and $k$ is an integer, $D(|A|^k)$ is the set of
$2k$-times  differentiable functions with square integrable $(2k)^{th}$
derivative. 
\end{remark}

The notion of weakly-coupled systems is closely related to the growth of the
expected value of the energy.
Here $\Re$
denotes the real part of a complex number.

\begin{definition}\label{DEF_weakly_coupled}
Let $k$ be positive  and let  $(A,B_1,\ldots,B_p)$ satisfy
Assumption~\ref{ass:ass}.\ref{EgaliteDomaine}. Then $(A,B_1,\ldots,B_p)$ is
\emph{$k$-weakly-coupled}
if for every $u\in {\mathbf R}^p$, $D(|A+\sum_l u_lB_l|^{k/2})=D(|A|^{k/2})$ and
 there exists
a constant $C$ such that, for every
$1\leq l\leq p$, for every $\psi$ in $D(|A|^k)$, $ |\Re \langle |A|^k
\psi,B_l\psi \rangle |\leq C |\langle |A|^k \psi, \psi \rangle|$.

The \emph{coupling constant}  $c_k(A,B_1,\ldots,B_p)$ of system $(A,B_1,\ldots,B_p)$ of
order $k$ is the quantity
$$
\sup_{\psi\in D(|A|^{k})\setminus\{0\}} \sup_{1\leq l \leq p} \frac{ |\Re
\langle |A|^k
\psi,B_l\psi \rangle |}{|\langle |A|^k \psi, \psi \rangle|}.
$$
\end{definition}
\begin{remark}
The terminology \emph{weak-coupling} refers to the weakness of  $B$ in the scale of $A$. In other words the effect of $B$ on the spectral properties of $A$ is small enough to have a weak coupling effect on the
Galerkin approximations associated with eigenvectors of $A$ (see Lemma~\ref{PRO_troncature} below) or the  boundedness in the
$s$-norm of $A$ of the evolution (see Proposition~\ref{PRO_croissance_norme_A} below).  The weakness of this action can also be seen
through the transition probabilities or energy transitions between eigenstates (see Lemma~\ref{PRO_fake_Taylor} below).
\end{remark}

We have the following technical interpolation result proved in Appendix~\ref{appA}.
\begin{lemma}\label{lem:Interpolation}
Let $A$ and $A'$ be invertible (from their respective domains to
$H$) skew-adjoint operators with compact
resolvent. Let $k$ be a positive real. Assume that
$D(|A|^k)=D(|A'|^k)$. Then
for any real $s\in(0,k)$, $D(|A|^s)=D(|A'|^s)$.
\end{lemma}

A first property of the propagator of a weakly-coupled system is given by the following proposition whose proof is  in  Appendix~\ref{appB}.

\begin{proposition}\label{PRO_croissance_norme_A} Let $k$ be a positive number
and let $(A,B_1,\ldots,B_p)$ satisfy Assumption~\ref{ass:ass} and be
$k$-weakly-coupled.  Then,
for every $\psi_{0} \in D(|A|^{k/2})$, $K>0$,
 $T\geq 0$, and piecewise constant function $u=(u_1,\ldots,u_p)$  for which
$\sum_{l=1}^p\|u_l\|_{L^1} \leq  K$, one has
\begin{equation}
\left\|\Upsilon^{u}_{T}(\psi_{0})\right\|_{k/2} \leq 
e^{c_{k}(A,B_1,\ldots,B_p) K} \| \psi_0 \|_{k/2}.
\end{equation}
\end{proposition}

\subsection{Good Galerkin approximation}

In this section we show that a weakly-coupled system admits a finite
dimensional approximation  with  trajectories  close, at
any time, to the solutions of the original
infinite dimensional system.
For every $N$ in $\mathbf{N}$, we define the orthogonal projection
$$
 \pi_N:\psi \in H \mapsto \sum_{j\leq N} \langle \phi_j,\psi\rangle
\phi_j \in H.
$$
\begin{lemma}\label{PRO_troncature}
Let $k$ be a positive number, $(A,B_1,\ldots,B_p)$ satisfy
Assumption~\ref{ass:ass}, and be $k$-weakly-coupled. Assume that there exist $d
>0$, $0 \leq r <k$ such that $\|B_l\psi \|\leq d \| \psi \|_{r/2}$ for every
$\psi$ in $D(|A|^{r/2})$ and $l$ in $\{1,\ldots,p\}$. Then,
for every  $K\geq 0, n\in \mathbf{N}$, $N \in \mathbf{N}$,
$(\psi_j)_{1\leq j \leq n}$ in $D(|A|^{k/2})^n$,
and for every piecewise constant function $u=(u_1,\ldots,u_p)$,
such that $\sum_{m=1}^p\|u_{m}\|_{L^1} \leq K$, one has
\begin{equation}\label{eq:feps2}
\|B_l(\mathrm{Id} - \pi_{N})
\Upsilon^{u}_{t}(\psi_{j})\| \leq d\lambda_{N+1}^{(r-k)/2} e^
{ c_ { k } (A ,B_1,\ldots,B_p)K
}\|\psi_j\|_{k/2},
\end{equation}
for every $t \geq 0$, $l=1,\ldots,p$ and $j=1,\ldots,n$.
\end{lemma}

\begin{IEEEproof}
Fix $j \in \{1,\ldots,n\}$. For every $N > 1$, one has
\begin{eqnarray}\label{eq:estimates}
{\left\|(\mathrm{Id} - \pi_{N})
\Upsilon^{u}_{t}(\psi_{j})\right\|_{r/2}^2}&=&\sum_{n = N+1}^{\infty} \lambda_{n}^{r} | \langle  \phi_{n},
\Upsilon^{u}_{t}(\psi_{j}) \rangle|^{2}  \nonumber\\
&
 \leq&
\lambda_{N+1}^{r-k}\left\|\Upsilon^{u}_{t}(\psi_{j})\right\|_{k/2}^2.
\end{eqnarray}
By Proposition~\ref{PRO_croissance_norme_A},
$\left\|\Upsilon^{u}_{t}(\psi_{j})\right\|_{k/2}^2 \leq
e^{2c_{k}(A,B_1,\ldots,B_p)K}\| \psi_j
\|_{k/2}^2$ for
every $t>0$ and $u$  of $L^1$-norm smaller than $K$. Equation
\eqref{eq:feps2} follows as, for every $l=1,\ldots,p$,
$\|B_l\psi \|\leq d \| |A|^{\frac{r}{2}}\psi \|$.\marginpar{\tiny }
\end{IEEEproof}

\begin{remark}
 Since $r<k$, then $\left\|B_l(\mathrm{Id} - \pi_{N})
\Upsilon^{u}_{t}(\psi_{j})\right\|_{r/2}$
tends to $0$,  uniformly with
respect to $u$, as $N$ tends to infinity.
\end{remark}

\begin{definition}
Let $N \in \mathbf{N}$.  The \emph{Galerkin approximation}  of \eqref{eq:main}
of order $N$ is the system in $H$
\begin{equation}\label{eq:sigma}
\dot x = \left(A^{(N)} + \sum_{l=1}^{p}u_{l} B_{l}^{(N)}\right) x \tag{$\Sigma_{N}$}
\end{equation}
where $A^{(N)}=\pi_N A \pi_N$ and $B_{l}^{(N)}=\pi_N B_{l} \pi_N$ are the
\emph{compressions} of $A$ and $B_{l}$ (respectively).
\end{definition}

We denote by $X^{u}_{(N)}(t,s)$ the propagator of \eqref{eq:sigma}
for a $p$-uple of piecewise constant functions $u = (u_{1}, \ldots, u_{p})$.

\begin{remark}
The operators $A^{(N)}$ and $B_l^{(N)}$ are defined
on the infinite dimensional space $H$. However, they have finite rank and
the dynamics of $(\Sigma_N)$ leaves invariant the $N$-dimensional space
$\mathcal{L}_{N} = \mathrm{span}_{1\leq j \leq N} \{\phi_j\}$. Thus,
$(\Sigma_N)$ can be seen as a finite dimensional bilinear dynamical system in
$\mathcal{L}_{N}.$
\end{remark}

 The operator $A$ written in one of its eigenvector basis is diagonal and its
dynamics is decoupled on each eigenspace. Thus the
projection of the dynamics coincides with
the dynamics of the associated truncation. 

One of the most important consequence of the weak-coupling assumption is that, even though the coupling action of the operator $B$ can give rise to intricate dynamics, this action is weak enough to allow approximations by the
dynamics of the truncations as stated in the theorem below.
\begin{thm}[Good Galerkin Approximation]\label{PRO_Good_Galerkin_approximation}
Let $k$ and $s$ be non-negative  numbers with
$0\leq s <k$. Let $(A,B_1,\ldots,B_p)$
satisfy Assumption~\ref{ass:ass} and be $k$-weakly-coupled. Assume that there
exist $d>0$ and $0\leq r<k$ such that $\|B_l\psi \|\leq d \|\psi \|_{r/2}$ for
every $\psi$ in $D(|A|^{r/2})$ and $l$ in $\{1,\ldots,p\}$. Then
for every $\varepsilon > 0 $, $K\geq 0$, $n\in \mathbf{N}$, and
$(\psi_j)_{1\leq j \leq n}$ in $D(|A|^{k/2})^n$
there exists $N \in \mathbf{N}$
such that
for every piecewise constant function $u=(u_1,\ldots,u_p)$
$$
\sum_{l=1}^p\|u_l\|_{L^{1}} < K \implies \| \Upsilon^{u}_{t}(\psi_{j}) -
X^{u}_{(N)}(t,0)\pi_{N} \psi_{j}\|_{s/2} < \varepsilon,
$$
for every $t \geq 0$ and $j=1,\ldots,n$.
\end{thm}

\begin{IEEEproof}
Consider the case $s=0$.
Fix $j$ in $\{1,\ldots,n\}$ and consider
the map $t\mapsto
\pi_{N} \Upsilon^{u}_{t}(\psi_{j})$ that is absolutely continuous and satisfies,
for almost every $t \geq 0$,
\begin{eqnarray*}
\frac{\mathrm{d}}{\mathrm{d}t} \pi_{N} \Upsilon^{u}_{t}(\psi_{j}) &= &(A^{(N)} + \sum_{l=1}^p u_l
B_l ^{(N)}) \pi_{N} \Upsilon^{u}_{t}(\psi_{j}) \\
&& \quad + \sum_{l=1}^p
u_l (t) \pi_{N} B_l (\mathrm{Id} - \pi_{N}) \Upsilon^{u}_{t}(\psi_{j}).
\end{eqnarray*}
Hence, by variation of constants, for every $t \geq 0$,
\begin{eqnarray}\label{EQ_preuve_good_Galerkin}
\lefteqn{\pi_{N} \Upsilon^{u}_{t}(\psi_{j})= X^{u}_{(N)}\!(t,0)  \pi_{N}\psi_{j}}  \nonumber \\
&& +
\sum_{l=1}^p \! \!\int_{0}^{t} \!\!\!\!
X^{u}_{(N)}\!(t,s) \pi_{N} B_l(\mathrm{Id} - \pi_{N}) \Upsilon^{u}_{s}(\psi_{j})
u_l(\tau)  \mathrm{d}\tau.
\end{eqnarray}
By Lemma~\ref{PRO_troncature}, the norm of $t \mapsto B_l(\mathrm{Id} -
\pi_{N}) \Upsilon^{u}_{t}(\psi_{j})$ is
less than
$d\lambda_{N+1}^{(r-k)/2}e^{c_{k}(A,B_1,\ldots,B_p)K}\|\psi_{j}\|_
{k/2}$. Since $X^{u}_{(N)}(t,s) $ is unitary,
\begin{eqnarray}\label{eq:remark}
\lefteqn{\|\pi_{N} \Upsilon^{u}_{t}(\psi_{j}) -X^{u}_{(N)}(t,0) \pi_{N}\psi_{j}\|} \nonumber \\
&& \leq
K
d\lambda_{N+1}^{(r-k)/2}e^{c_{k}(A,B_1,\ldots,B_p)K}\|\psi_{j}\|_{k/2}.
\end{eqnarray}
Then
\begin{eqnarray}
\lefteqn{\|\Upsilon^u_t(\psi_j) - X^u_{(N)}(t,0)\pi_{N} \psi_{j}\|} \nonumber\\
& \leq& \!\! \|(\mathrm{Id} - \pi_{N}) \Upsilon^{u}_{t}(\psi_{j})\|
\!+\! \|\pi_{N} \Upsilon^{u}_{t}(\psi_{j}) -\! X^{u}_{(N)}(t,0)  \pi_{N} \psi_{j}
\| \nonumber \\
& \leq& \lambda_{N+1}^{-k/2}e^{c_{k}(A,B_1,\ldots,B_p)K}\|\psi_{j}\|_{k/2} \nonumber \\
&& \quad 
+Kd\lambda_{N+1}^{(r-k)/2}e^{c_{k}(A,B_1,\ldots,B_p)K}\|\psi_{j}\|_{k/2}.
\label{eq:qwer}
\end{eqnarray}
This completes the proof for $s=0$  since $\lambda_N$ tends to infinity as $N$
goes to infinity.

Note that, if $\mathcal X$ is a set and $(v_n)_{n \in \mathbf{N}}$ is a sequence
of functions from ${\mathcal
X}$ to $H$ that tends uniformly to $0$ (the null function)
for the $s_1$-norm and it is uniformly bounded for the $s_2$-norm for $s_1<s_2$, then
$(v_n)_{n \in \mathbf{N}}$ tends uniformly to $0$ in
the ${\frac{s_1+s_2}{2}}$-norm. This is a consequence of Cauchy--Schwarz
inequality, indeed
\begin{eqnarray}\label{eq:interpolation}
\|v_n \|^2_{\frac{s_1+s_2}{2}}&=&|\langle |A|^{\frac{s_1+s_2}{2}} v_n,
|A|^{\frac{s_1+s_2}{2}}  v_n \rangle | \nonumber \\
&=&|\langle |A|^{s_1} u_n, |A|^{s_2}  v_n \rangle
|\leq\|v_n\|_{s_1} \|v_n\|_{s_2}.
\end{eqnarray}
To conclude the proof in the
general case $0<s<k$, we apply iteratively this interpolation result with
$v_N: (t,u) \mapsto (X^u_{(N)}(t,0)\pi_{N}-\Upsilon^u_t)\psi_j$, defined on
${\mathcal X}=[0,+\infty) \times \{u \in L^1\,:\, \|u\|_{L^1} \leq
K\}$.  From the first part of the proof, $(v_N)_N$ tends uniformly to zero for the $s_1=0$ norm and it is bounded for the $s_2=k$ norm. Hence {by~\eqref{eq:interpolation}}, the sequence $(v_N)_N$ tends uniformly to zero for the $k/2$ norm. Applying once again the interpolation estimate~{\eqref{eq:interpolation}} with $s_1=k/2$ and $s_2=k$, we obtain that the sequence $(v_N)_N$ tends uniformly to zero for the $3k/4$ norm. After $l$ interpolations,   we obtain that the sequence $(v_N)_N$ tends uniformly to zero for the $k(1-1/2^l)$ norm. Conclusion follows by choosing an integer $l$ such that  $k(1-1/2^l)>s$. 
\end{IEEEproof}

\begin{remark}\label{rk:estimates}
In the case $s=0$, there is an explicit estimate for the order of the Galerkin
approximation  which existence is
stated in Theorem~\ref{PRO_Good_Galerkin_approximation}. For instance,
by~\eqref{eq:remark}, $\| \pi_N \Upsilon^{u}_{t}(\psi_{j}) -
X^{u}_{(N)}(t,0)\pi_{N} \psi_{j}\| < \varepsilon$ if $N$ is such that, for $j=1,\ldots,n$,
\begin{equation}\label{EQ_estimation_N_projection}
\lambda_{N+1} > \left(
\frac{Kde^{c_{k}(A,B_1,\ldots,B_p)K}\|\psi_{j}\|_{k/2}}{\varepsilon} \right
)^{\frac{2}{k-r}} \!\!\!\!\!\!\!.
\end{equation}
\end{remark}

\subsection{Approximate controllability in $s$-norm}
Let $(A,B_1,\ldots,B_p)$ be a $k$-weakly-coupled system. For every $\phi$ in $ D(|A|^{k/2})$, every $T\geq 0$ and every piecewise constant function $u:[0,T]\rightarrow \mathbf{R}^p$, 
one has $\Upsilon^u_T \phi \in D(|A|^{k/2})$, which is a deep   obstruction  to exact
controllability. But this property also  provides powerful tools for the study of the approximate controllability.

\begin{definition}
 Let $(A,B)$ satisfy Assumption~\ref{ass:ass}.
 A subset $S$ of
$\mathbf{N}^2$ \emph{couples} two levels $j,l$ in $\mathbf{N}$,
if
there exists a finite sequence $\big
((s^{1}_{1},s^{1}_{2}),\ldots,(s^{q}_{1},s^{q}_{2}) \big )$
in $S$ such that
\begin{description}
\item[$(i)$] $s^{1}_{1}=j$ and $s^{q}_{2}=l$;
\item[$(ii)$] $s^{j}_{2}=s^{j+1}_{1}$ for every $1 \leq j \leq q-1$;
\item[$(iii)$] $\langle  \phi_{s^{j}_{1}}, B \phi_{s^{j}_{2}}\rangle \neq 0$
for $1\leq j \leq
q$.
\end{description}

The subset $S$ is called a \emph{connectedness chain}   for $(A,B)$ if $S$
couples every pair of levels in $\mathbf{N}$.
A connectedness chain is said to be \emph{non-degenerate} (or sometimes \emph{non-resonant}) if for every $(s_1,s_2)$
in $S$, $|\lambda_{s_1}-\lambda_{s_2}|\neq |\lambda_{t_1}-\lambda_{t_2}|$ for
every $(t_{1},t_{2})$ in $\mathbf{N}^2\setminus\{(s_1,s_2),(s_2,s_1)\}$ such
that $\langle \phi_{t_{2}}, B \phi_{t_{1}}\rangle  \neq 0$.
\end{definition}

\begin{definition}
Let $(A,B)$ satisfy Assumption~\ref{ass:ass} and $s>0$. The system $(A,B)$ is
approximately simultaneously controllable for the $s$-norm if for every
$\psi_1,\ldots,\psi_n\in D(|A|^{s})$, $\hat{\Upsilon} \in U(H)$ such that $\hat\Upsilon(\psi_{1}),\ldots,\hat\Upsilon(\psi_{n}) \in D(|A|^{s})$, 
and
$\varepsilon>0$, there  exists a piecewise constant function
$u_\varepsilon: [0,T_\varepsilon]\to \mathbf{R}$ such that
$$
\|\hat{\Upsilon}\psi_j-\Upsilon_{T_\varepsilon}^{u_\varepsilon}
\psi_j\|_s<\varepsilon.
$$
for every $j=1,\ldots,n$.
\end{definition}

\begin{proposition}\label{prop:controlHs}
Let $k$ be a positive number. Let $(A,B)$  satisfy Assumption~\ref{ass:ass}, be
$k$-weakly-coupled, and admit a non-degenerate chain of connectedness. Assume that
there exist $d>0$, $0\leq r<k$ such that $\|B \psi \|\leq d \|
|A|^{\frac{r}{2}}\psi \|$, for every $\psi$ in $D(|A|^{\frac{r}{2}})$.
Then $(A,B)$ is approximately simultaneously controllable
for the norm $\|\cdot\|_{{s/2}}$ for every $s<k$.
\end{proposition}

\begin{IEEEproof}
Fix $\varepsilon >0$, $\psi_{1},\ldots,\psi_{n} \in D(|A|^{s/2})$, and
$\hat{\Upsilon} \in U(H)$ such that
$\hat{\Upsilon}(\psi_{1}),\ldots,\hat{\Upsilon}(\psi_{n})\in D(|A|^{s/2})$.
Fix $n_{1}$ sufficiently large such that $\|\hat{\Upsilon}(\psi_{j})
- \pi_{n_{1}}\hat{\Upsilon}(\psi_{j})\|_{s/2} < \varepsilon/ 3$ for
every $j=1,\ldots,n$.

There exist $l_{1},\ldots,l_{n}$ such that $t \mapsto (e^{it\lambda_{l_{1}}},
\ldots,e^{it\lambda_{l_{n}}})$ is $\varepsilon$-dense in the torus
$\mathbf{T}^{n}$
(see~\cite[Proposition~6.1]{Schrod2}). Call
$m=\max\{n_{1},l_{1},\ldots,l_{n}\}$.

By~\cite{Schrod2} the existence of a non-degenerat chain of connectedness is sufficient 
for the approximate controllability of~$(A,B)$ in the norm of $H$. More precisely,
by~\cite[Remark 5.9]{Schrod2} there exists $K_{1}$ such that for every $\eta >
0$ there exist a control $u_{1}^{\eta}$  satisfying $\|u^{\eta}_{1}\|_{L^{1}}
\leq K_{1}$ and $\theta_{1},\theta_{2}, \ldots, \theta_{n}$, such that
$\|\Upsilon^{u^{\eta}_{1}}_{T_{1}}(\psi_{j}) - e^{\mathrm{i} \theta_{j}}
\phi_{l_{j}}\| < \eta, $ for every $j=1,\ldots,n$.

Similarly, since as shown in~\cite[Section 6.1]{Schrod2} the hypotheses sufficient  for controllability
(and in particular the one of~\cite[Remark 5.9]{Schrod2}) hold for the
system $(-A,-B)$, we have existence
of $K_{2}$ such that for every $\eta > 0$ there exists
$u_{2}^{\eta}$
satisfying $\|u^{\eta}_{2}\|_{L^{1}} \leq K_{2}$ and  $\bar \theta_{1} ,\ldots,
\bar \theta_{n} \in\mathbf{R}$ such that
the solution of the system
$$
\frac{\mathrm{d}\psi}{\mathrm{d}t}(t) = - (A+u(t) B)\psi(t)
$$
at time $T_{2}$ with initial state $\hat{\Upsilon}(\psi_{j})$
and
corresponding to the control $u_{2}^{\eta}$ is $\eta$-close in the norm of $H$
to
$e^{\mathrm{i} \bar\theta_{j}} \phi_{l_{j}}$ for every $j=1,\ldots,n$.

Let $\tau$ such that $\|e^{\mathrm{i} \tau \lambda_{l_{j}}} e^{i
\theta_{j}} - e^{\mathrm{i} \bar\theta_{j}}\|<\eta$ for every $j=1,\ldots,n$.
Let $T = T_{1}+\tau+T_{2}$ and let  $u:[0,T] \to \mathbf{R}$ be the piecewise
constant control defined by
$$
u^\eta(t)=
\left\{
\begin{array}{ll}
u^{\eta}_{1}(t) & t \in [0,T_{1}),\\
0 & t \in [T_{1},T_{1}+\tau),\\
u^{\eta}_{2}(T_{2} - (t - T_{1} - \tau)) & t \in [T_{1}+\tau,T],
\end{array}
\right.
$$
The control $u^\eta$ above steers a solution of $\dot{\psi} = (A+uB)\psi$ with
initial state $\psi_{j}$, $3\eta$-close in the norm $\|\cdot\|$
to
$\hat{\Upsilon}(\psi_{j})$ in a time $T$, namely
$
\| \hat{\Upsilon}(\psi_{j}) -
\Upsilon^{u^{\eta}}_{T}(\psi_{j}) \|\leq   3 \eta.
$

Let $K= K_{1} + K_{2}$. By Lemma~\ref{PRO_troncature}, we have that there
exists $N=N(\varepsilon,K,s) > n$ such that
$$
\|u\|_{L^{1}} \leq K \implies \|(\mathrm{Id} -
\pi_{N})\Upsilon_t^{u}(\psi_{j})\|_{{s/2}} < \frac \varepsilon 3\,,
$$
for every $j=1,\ldots,n$ and $t \geq 0$.

Note that, on $\mathrm{span}\{\phi_{1},\ldots,\phi_{N}\}$, we have
$\|\cdot\|_{s/2} \leq\lambda^{s/2}_{N} \|\cdot\|$.
Therefore for every $j =1, \ldots ,n$,
\begin{eqnarray*}
\lefteqn{\|\hat{\Upsilon}(\psi_j)  - \Upsilon^{u^{\eta}}_{T}(\psi_{j})
\|_{s/2}}\\
 & \leq&     \| (\mathrm{Id} - \pi_{N})(\hat{\Upsilon}(\psi_{j}) -
\Upsilon^{u^{\eta}}_{T}(\psi_{j}) )  \|_{{s/2}}  \\
&&+\|
\pi_{N}(\hat{\Upsilon}(\psi_{j}) - \Upsilon^{u^{\eta}}_{T}(\psi_{j}) )
\|_{{s/2}}\\
&\leq&   \| (\mathrm{Id} - \pi_{N})\hat{\Upsilon}(\psi_{j})\|_{{s/2}} + \|
(\mathrm{Id} - \pi_{N})\Upsilon^{u^{\eta}}_{T}(\psi_{j}) )\|_{{s/2}}\\
&& +
\lambda^{s/2}_{N} \| \hat{\Upsilon}(\psi_{j}) -
\Upsilon^{u^{\eta}}_{T}(\psi_{j}) \|\\
&\leq &  \frac{2\varepsilon}{3} + 3 \lambda^{s/2}_{N}\eta < \varepsilon,
\end{eqnarray*}
for $\eta$ sufficiently small.
\end{IEEEproof}

\section{The bounded case}\label{SEC_compact_case}

\begin{proposition}\label{prop:wc123}
Let $k$ be a positive integer. Assume that for every $u\in {\mathbf R}^p$,
$D(|A|^{\frac{k}{2}})=D(|A+\sum_lu_lB_l|^{\frac{k}{2}})$ and that  for every
$l=1,\ldots,p$ the restriction of $B_l$ to $D(|A|^{\frac{k}{2}})$ is bounded
for the ${\frac{k}{2}}$-norm. Then $(A,B_1,\ldots,B_p)$ is
$k$-weakly-coupled.
\end{proposition}
\begin{IEEEproof} For every $l=1,\ldots,p$, let $\|B_l \psi\|_{k/2} \leq
C_{l,k}\|\psi\|_{k/2}$ for every $\psi$ in
$D(|A|^{k})$.
Then $|\langle A^k \psi, B_l \psi \rangle|=|\langle  |A|^{\frac{k}{2}} \psi,
|A|^{\frac{k}{2}} B_l \psi \rangle| \leq \||A|^{\frac{k}{2}} \psi \|
\||A|^{\frac{k}{2}} B_l \psi\| \leq C_{l,k} \||A|^{\frac{k}{2}}\psi \|^2
=C_{l,k} |\langle A^k \psi,  \psi \rangle|$ for every $\psi$ in
$D(|A|^{k})$.\end{IEEEproof}

\subsection{Example: single trapped ion}\label{SEC_trap_ion}

This example is a model  of a two level ion trapped in a harmonic potential 
and under the action of an external field. This model has been extensively studied (see for example~\cite{eberly-law},  \cite{Boscain_Adami},  \cite{puel}, and \cite{rangan}).

The state of the system is $(\psi_e,\psi_g)$ in
$H=L^2(\mathbf{R},\mathbf{C})\times L^2(\mathbf{R},\mathbf{C})$. The dynamics
is
given by {two coupled harmonic oscillators}
$$
\left \{ \begin{array}{lcl}
\mathrm{i}\frac{\partial \psi_e}{\partial t}&=& \omega (-\Delta +x^2)\psi_e +
\Omega \psi_e\\
&& + \left( u_1(t)\cos({\sqrt{2}\eta x})+  u_2(t)\sin({\sqrt{2}\eta
x})\right) \psi_g\\
\mathrm{i}\frac{\partial \psi_g}{\partial t}&=& \omega (-\Delta +x^2)\psi_g +
\Omega \psi_g \\&&+\left( u_1(t)\cos({\sqrt{2}\eta x})+  u_2(t)\sin({\sqrt{2}\eta
x})\right) \psi_e
\end{array}
\right.
$$
where $\omega, \Omega, \eta$ are positive constants related to the
physical properties of the system. The two real valued controls $u_1$ and $u_2$
are usually a sum of periodic functions with positive frequencies $\Omega$,
$\Omega+\omega$ and $\Omega-\omega$. With our notations, the dynamics reads
\begin{equation}\label{eq:1231}
\frac{\mathrm{d}\psi}{\mathrm{d}t}=A \psi +u_1(t) B_1\psi + u_2(t) B_2\psi
\end{equation}
where $A$ is the diagonal operator $A:(\psi_e,\psi_g)\mapsto -{\mathrm i
}(\omega (-\Delta +x^2)\psi_e+\Omega \psi_e, \omega (-\Delta +x^2)\psi_g+
\Omega\psi_g)$,  $B_1:(\psi_e,\psi_g)\mapsto -{\mathrm i
}(\cos({\sqrt{2}\eta x}) \psi_g,\cos({\sqrt{2}\eta
x}) \psi_e)$, and $B_2: (\psi_e,\psi_g)\mapsto -{\mathrm i
}( \sin ({\sqrt{2}\eta x}) \psi_g, \sin ({\sqrt{2}\eta x}) \psi_e)$.

By~\cite[Theorem XIII.69 and Theorem XIII.70]{reed-simon-4}, the operator $A$ is
skew-adjoint with discrete spectrum and admits a family of eigenfunctions which
forms an orthonormal basis of $H$. Since $B_1$ and $B_2$ are bounded then, for
every real  constants $u_1$ and $u_2$, $A+u_1 B_1 + u_2 B_2$ is
skew-adjoint with the same domain of $A$ (see \cite[Theorem
X.12]{reed-simon-2}).
The spectrum of $A$ is the sequence $(-\mathrm{i} \lambda_n)_{n \in
\mathbf{N}}=-\mathrm{i} (\omega(n+1/2) + \Omega)_{n\in
\mathbf{N}}$. For every $n$ in $\mathbf{N}$, the eigenvalue
$-\mathrm{i}\lambda_n$ has multiplicity 2 and is associated with the
2-dimensional subspace of $L^2(\mathbf{R},\mathbf{C})\times L^2(\mathbf{R},\mathbf{C})$ spanned by $\{ (f_n,0),
(0,f_n)\}$ where $f_n$ is the $n^{th}$ Hermite function.
Assumption~\ref{ass:ass}
is then verified.
Since, for every $k$ in $\mathbf{N}$, all derivatives up to order $k$ of $x
\mapsto \cos({\sqrt{2}\eta x})$ and $x \mapsto \sin({\sqrt{2}\eta x})$ are
bounded for the $L^\infty$-norm  by $C_k=2^{\frac{k}{2}} \eta^k$ on $\mathbf{R}$
then  $B_1$ and $B_2$ are bounded
by ${2^k}C_k$
on
$D(|A|^\frac{k}{2})$ for every $k$.
Moreover for every
$(u_1,u_2)\in{\mathbf R}^2$, $D(A^k)=D((A+u_1B_1+u_2B_2)^k)$. Indeed by induction
on $k$
\begin{eqnarray*}
\lefteqn{D((A+u_1B_1+u_{2}B_2)^{k+1})}\\
&=&\{\psi\in D((A+u_1B_1+u_2B_2)^k)\,:\, \\
&&(A+u_1B_1+u_2B_2)\psi\in D((A+u_1B_1+u_2B_2)^k)\}\\
&=&\{\psi\in D(A^k)\,:\, (A+u_1B_1+u_2B_2)\psi\in D(A^k) \}\\
&=& D(A^{k+1}),
\end{eqnarray*}
since $(u_1B_1+u_{2}B_2)\psi \in D(A^k)$ when $\psi \in D(A^k)$.
Hence for every
$(u_1,u_2)\in{\mathbf R}^2$, $D(|A|^k)=D(|A+u_1B_1+u_2B_2|^k)$ and
Lemma~\ref{lem:Interpolation} provides $D(|A|^s)=D(|A+u_1B_1+u_2B_2|^s)$ for any
$s>0$. Hence, by Proposition~\ref{prop:wc123} the system $(A,B_1,B_2)$ is
$k$-weakly-coupled for every $k$, with coupling constant smaller than $2^k
C_k$.

\subsection{The case of a compact manifold}\label{sec:compact}
We focus here on the case where the space $\Omega$ is a compact Riemannian
manifold (without boundary).  By Rellich-Kondrakov and
Weyl theorems, if $V$ is essentially bounded the operator $A=-\mathrm{i}(\Delta
+V):H^2(\Omega)\rightarrow L^2(\Omega,\mathbf{C})$ has purely discrete spectrum
 $(-{\mathrm i }\lambda_n)_{n\in\mathbf{N}}$ with $\lambda_n$ non-decreasing to
infinity (see for instance \cite[Theorem 7.2.6]{buser}). Note that $\lambda_1$
is not necessarily positive but this
can be assumed considering $A+\mathrm{i}(\lambda_1-1)$  instead of $A$.
This shift gives a physically irrelevant phase term,
$e^{\mathrm{i}t(\lambda_1-1)}$, on the dynamics associated with $A$.

\begin{lemma}\label{lem:domainA}
Let $k$ be a positive integer, $\Omega$ be a compact Riemannian manifold,
$V:\Omega\rightarrow \mathbf{R}$ be  $C^{2 k}(\Omega)$.  Then the domain of the
operator
 $(\Delta +V)^k$
is $H^{2k}(\Omega)$.
\end{lemma}

\begin{IEEEproof}
Since $\Omega$ is compact it is sufficient to prove the proposition on a
bounded domain of $\mathbf{R}^n$.
The operator $-\mathrm{i}A= \Delta +V$ is an elliptic operator of order $2$. By
\cite[Theorem 8.10]{gilbarg} if $Af \in H^k(\Omega)$ then $f \in
H^{k+2}(\Omega)$ and by
induction we have that $D(|A|^k) = H^{2 k}(\Omega)$.
\end{IEEEproof}

\begin{proposition}\label{pro_weak_coupled_compact}
Let $k$ be a positive integer, $\Omega$ be a
compact Riemannian manifold,
$V,W:\Omega\rightarrow \mathbf{R}$ be two $C^{2 k}(\Omega,\mathbf{R})$ functions
on $\Omega$.
Define $A=-\mathrm{i}(\Delta +V):D(A)\rightarrow L^2(\Omega,\mathbf{C})$ and
$B=\mathrm{i}W:L^2(\Omega,\mathbf{C})\rightarrow L^2(\Omega,\mathbf{C})$.
Then $(A,B)$ is $k$-weakly-coupled.
\end{proposition}

\begin{IEEEproof}
Note that for every $f \in C^{2k}$ there exists a constant $C_k=2^{2k+1}
\sup_{0\leq j\leq 2k} \|W^{(j)}\|_{L^{\infty}(\Omega,\mathbf{R})}$ such that
$\|W f\|_{H^{2k}} \leq C_k \|f\|_{H^{2k}}$. From Lemma~\ref{lem:domainA}, the
norm $\|\cdot\|_{H^{2k}}$ and the $k$-norm are equivalent.
Therefore, by Proposition~\ref{prop:wc123}, the
system is $k$-weakly-coupled.
\end{IEEEproof}

\begin{remark}
As a consequence of Lemma~\ref{lem:domainA}
and Proposition~\ref{pro_weak_coupled_compact} we have that, in the case
of a
compact manifold, if the potentials are in $C^{m}(\Omega)$ then
Theorem~\ref{PRO_Good_Galerkin_approximation} applies with $k=m/2-1$ and
$r=0$.
\end{remark}

\subsection{Example: orientation of a rotating molecule in the plane}
\label{sec:molecule2d}

We consider  a rigid bipolar molecule rotating in a plane. Its only
degree of freedom is the rotation around its centre of mass.  The molecule
is submitted to an electric field of constant direction with variable intensity
$u$.  The orientation of the molecule is an angle in $\Omega=SO(2) \simeq
\mathbf{R}/2\pi \mathbf{Z}$. The dynamics is governed by the Schr\"odinger
equation
$$
\mathrm{i}\frac{\partial \psi(\theta,t) }{\partial t} =
\left(-\frac{\partial^2}{\partial \theta^2} + u(t)\cos \theta \right)
\psi(\theta, t), \quad \theta \in \Omega.
$$
Note that the parity (if any) of the wave function is preserved by the above
equation.  We consider then the Hilbert space $H=\{\psi\in
L^2(\Omega,\mathbf{C}): \psi \mbox{ odd } \}$, endowed with the Hilbert product
$\langle f,g\rangle =\int_{\Omega} \bar{f}g$. The eigenvalue of the skew-adjoint
operator $A = \mathrm{i}\frac{\partial^2}{\partial \theta^2}$
associated with the eigenfunction $\phi_k:\theta\mapsto\sin(k\theta)/\sqrt{\pi}$
is $-\mathrm{i} \lambda_k = - \mathrm{i} k^2$, $k \in \mathbf{N}$.
The domain of $|A|^{k}$ is the Hilbert space
$H^k_e=\{\psi\in H^{2k}(\Omega,\mathbf{C}): \psi \mbox{ odd }
\}$. The skew-symmetric operator $B=-\mathrm{i} \cos \theta$ is bounded on
$D(|A|^{k/2})$ for every $k$.
By Proposition~\ref{prop:wc123},  for every $k$ in $\mathbf{N}$, $(A,B)$ is
$k$-weakly-coupled.
Theorem~\ref{PRO_Good_Galerkin_approximation} applies for every $k$ with
 $r=0$ and $d= 1$.
 In Section \ref{SEC_Rotation_Tridiagonal} we also give an estimate on the
 coupling constant $c_k(A,B)$  for this system.

From the viewpoint of the controllability problem, notice that the operator $B$ couples only
adjacent eigenstates, that is $\langle
\phi_l, B \phi_j \rangle = 0$ if and only if $|l-j|>1$. Since $\lambda_{l+1} -
\lambda_l = 2l+1$ then $\{(j,l)\in \mathbf{N}^2 \,:\, |l-j|=1\}$ is
a non-degenerate connectedness chain for $(A,B)$. Therefore, by Proposition
\ref{prop:controlHs} the system provides an example of approximately
simultaneously controllable system in norm $H^k(\Omega)$ for every
$k$. Note that, since the eigenstates belong to $H^{k}(\Omega)$ for every $k$ then
the
reachable set from any eigenstate is contained in $H^{k}(\Omega)$ for every
$k$.
\subsection{Example: orientation of a rotating molecule in the space}
We present the physical example of a rotating rigid bipolar molecule. Unlike
last example the motion of the molecule is not confined to a plane.  The model
then can be represented by the Schr\"odinger equation on the sphere.
In this case, $\Omega=\mathbf{S}^2$ is the unit sphere, the family
$(Y_{\ell}^m)_{\ell\geq 0, |m| \leq \ell}$ of the spherical harmonics is an
Hilbert basis of $H = L^2(\Omega,\mathbf{C})$, and the control is represented
by  three piecewise constant functions $u_1, u_2, u_3$. Using spherical coordinates $(\nu,\theta)$, the controlled
Schr\"odinger  equation is
\begin{eqnarray*}
\mathrm{i}\frac{\partial \psi(\nu, \theta,t)}{\partial t}&=
&\!\!-\Delta  \psi(\nu,
\theta,t)  + u_1(t)
\cos\theta\sin\nu ~\psi(\nu,
\theta,t) \\ && \!+\! \left ( u_2(t) \sin\theta\sin\nu +u_3(t) \cos \nu \right ) \psi(\nu,
\theta,t).
\end{eqnarray*}
Therefore, since $\Omega$ is compact,
Theorem~\ref{PRO_Good_Galerkin_approximation} applies for every integer $k$
with $d=1$ and $r=0$.

\section{Tri-diagonal systems}\label{SEC_tridiagonal}

We deal with the case where $p=1$ and $B$ couples
only adjacent levels of $A$.

\subsection{Tri-diagonal systems}

\begin{definition}
A system $(A,B)$ satisfying Assumption~\ref{ass:ass} is \emph{tri-diagonal} if
for every $j,k$ in $\mathbf{N}$,
$|j-k|>1$ implies $\langle \phi_j, B\phi_k \rangle =0$.
\end{definition}
{ In the following, we denote $b_{j,k}=\langle \phi_j, B\phi_k \rangle$.}

\begin{proposition}\label{PROP_Stability_Domain}
Assume that $(A,B)$ is tri-diagonal, that the sequence
$\left(\frac{\lambda_{n+1}}{\lambda_{n}}\right)_{n\in \mathbf{N}}$ is bounded,
and that the sequences
$\left (\frac{b_{n,n-1}}{\lambda_n}\right )_{n \in \mathbf{N}} , \left
(\frac{b_{n,n}}{\lambda_n} \right )_{n \in \mathbf{N}}$ tend to zero. Then, for
every $k$ in $\mathbf{N}$ and $u$ in $\mathbf{R}$, $D((A+uB)^k)=D(A^k)$.
Moreover, $D(A^k)$ is invariant for $e^{t(A+uB)}$ for any $u$ in $\mathbf{R}$ and
$t$ in $\mathbf{R}$.
\end{proposition}
\begin{IEEEproof}
The equality of $D((A+uB)^k)$ and $D(A^k)$ will follow from the Kato-Rellich
theorem (\cite[Theorem 1.4.2]{davies}). It suffices to check that for every $k$
in $\mathbf{N}$, $u$ in $\mathbf{R}$ and $\varepsilon>0$, there exists
$b_{\varepsilon}$ depending on $u$, $k$ and $\varepsilon$ such that, for every $\psi$ in
$D(A^k)$,  
\begin{equation}\label{EQ_recurrence_Ak}
\| ( (A+uB)^k -A^k) \psi \| \leq \varepsilon \|A^k \psi \| + b_{\varepsilon}
\|\psi\|.
\end{equation}
Let us prove that $B$ is bounded from $D(A^{r+1})$ to $D(A^r)$ for every
integer $r\geq 0$.
For every  $v$ in $D(A^r)$,
\begin{eqnarray*}
\|B v\|_{r}^2 \!\! &\!=\!& \!\sum_{n=1}^{\infty} \lambda_n^{2r} |\langle B v, \phi_n
\rangle | ^ 2
 =  \sum_{n=1}^{\infty} \lambda_n^{2r} |\langle  v, B\phi_n \rangle | ^ 2\\
& \!\leq\!&  \!\sum_{n=1}^{\infty} \bigg\{\lambda_n^{2r} ( |b_{n,n-1}|^2 |\langle 
\phi_{n-1}, v
\rangle | ^ 2
\!+\!|b_{n,n}|^2 |\langle  \phi_{n}, v \rangle | ^ 2 \\
&& \quad +  |b_{n,n+1}|^2 |\langle
\phi_{n+1}, v \rangle | ^ 2)\bigg\}
\end{eqnarray*}
\begin{eqnarray*}
~~~~~~& = &\sum_{n=1}^{\infty}\bigg\{
\lambda_{n-1}^{2r+2} \left(\frac{\lambda_{n}}{\lambda_{n-1}} \right)^{2r}
\frac{|b_{n,n-1}|^2}{\lambda_{n-1}^2}  |\langle  \phi_{n-1}, v \rangle | ^ 2\\
&& \quad 
+
 \lambda_{n}^{2r+2} \frac{|b_{n,n}|^2}{\lambda_{n}^2}  |\langle  \phi_{n}, v
\rangle | ^ 2\\
&& \quad 
+
\lambda_{n+1}^{2r+2} \left(\frac{\lambda_{n}}{\lambda_{n+1}} \right)^{2r}
\frac{|b_{n,n+1}|^2}{\lambda_{n+1}^2}  |\langle  \phi_{n+1}, v \rangle |^2\bigg\}.
\end{eqnarray*}
Now for every $\varepsilon >0$, let $n_0$ such that 
$\sup_{n \geq n_0}\frac{|b_{n,n}|^2}{\lambda_{n}^2} < \varepsilon/3$, 
$\sup_{n \geq n_0}\frac{|b_{n,n+1}|^2}{\lambda_{n+1}^2} < \frac{\varepsilon}{3}$, and $\sup_{n \geq
n_0}\frac{|b_{n,n-1}|^2}{\lambda_{n-1}^2}
<\frac{\varepsilon}{3C^{2r}}$ whith $C=\sup_n \lambda_{n+1}/\lambda_n$.
Note that the sequence $(\lambda_{n})_{n\in \mathbf{N}}$ is non-decreasing.
Then there exists $C_{\varepsilon}$ such that
\begin{eqnarray}\label{eq:BbddAk}
 \|B v\|_{r}^2   &\leq& \sum_{n=1}^{n_0}  \lambda_n^{2r} |\langle  v, B\phi_n
\rangle | ^ 2 +
\varepsilon \sum_{n\geq n_0-1} \lambda_{n}^{2r+2} |\langle  \phi_{n}, v \rangle
|^2 \nonumber \\
&\leq& C_\varepsilon \|v\|^2 +  \varepsilon \|v\|^{2}_{r+1}.
\end{eqnarray}

We prove \eqref{EQ_recurrence_Ak} by induction on $k$. For $k=1$ this is a
consequence of \eqref{eq:BbddAk} with $r=0$. The inductive step follows from the
fact that
\begin{eqnarray*}
(A+uB)^{k+1}-A^{k+1}&=&u((A+uB)^k B -A^k B) \\
&& +u A^kB  + ((A+uB)^k-A^k)A
\end{eqnarray*}
for
every $u$ in $\mathbf{R}$ and from~\eqref{eq:BbddAk}.
\end{IEEEproof}

\begin{proposition}\label{PRO_diagonal_s_weakly_coupled}
Let $(A,B)$ be a tri-diagonal system and let $k$ be a positive integer. Assume
that the sequence
$\left(\frac{\lambda_{n+1}}{\lambda_{n}}\right)_{n\in \mathbf{N}}$ is bounded,
that the sequences
$\left (\frac{b_{n,n-1}}{\lambda_n}\right )_{n \in \mathbf{N}} , \left
(\frac{b_{n,n}}{\lambda_n} \right )_{n \in \mathbf{N}}$ tend to zero, and that the
sequence $\displaystyle{\left(|b_{n,n+1}|\left (
\frac{\lambda_{n+1}^k}{\lambda_n^k}-1 \right )\right)_{n\in \mathbf{N}}}$ is
bounded. Then $(A,B)$ is $k$-weakly-coupled  with coupling constant
$$c_k(A,B)\leq \sup_n |b_{n,n+1}|\left (
\frac{\lambda_{n+1}^k}{\lambda_n^k}-1 \right ).$$
\end{proposition}

\begin{IEEEproof} 
For every $\psi$ in
$D(A)$, write $\psi=\sum_{j=1}^{\infty} x_j \phi_j$ where $x_j=\langle \phi_j,
\psi \rangle$. Since $\Re (b_{j,j}) =0$ then
\begin{eqnarray*}
\lefteqn{\Re \left( \langle |A|^k \psi, B \psi \rangle \right)}\\
  & =
& \Re \left( \sum_{j=1}^{\infty}
\lambda_{j}^k \bar{x}_{j}  b_{j+1,j}  x_{j+1} + \lambda_{j+1}^k \bar{x}_{j+1}
b_{j,j+1} x_j \right)\\
&=& \Re \left (\sum_{j=1}^{\infty} \lambda_j^k ( \bar{x}_j b_{j+1,j}  x_{j+1} -
x_j \bar{b}_{j+1,j} \bar{x}_{j+1} )  \right .\\
&& + (\lambda_{j+1}^k -\lambda_j^k)
\bar{x}_{j+1} b_{j,j+1} x_j \bigg)\\
&=& \Re \left( \sum_{j=1}^{\infty} (\lambda_{j+1}^k -\lambda_j^k) \bar{x}_{j+1}
b_{j,j+1} x_j \right)\\
&\leq &  \sum_{j=1}^{\infty} (\lambda_{j+1}^k -\lambda_j^k) |b_{j,j+1}|
\frac{|x_j|^2+|x_{j+1}|^2}{2}.
\end{eqnarray*}
By hypothesis, there exists $C$ such that
$|b_{j,j+1}|(\lambda_{j+1}^k-\lambda_j^k)\leq C \lambda_j^k$ for every $j$.
Hence, $|\Re \langle |A|^k \psi, B \psi \rangle| \leq   C \sum_{j=1}^{\infty}
\lambda_j^k |x_j|^2
 \leq C  \langle |A|^k \psi, \psi \rangle$.
The equality of the domains follows by Proposition~\ref{PROP_Stability_Domain}.
\end{IEEEproof}
{\color{red}
}

\subsection{Estimates for tri-diagonal systems}

{
\begin{lemma} \label{PRO_fake_Taylor}
Let $(A,B)$ be a tri-diagonal system and $n<l$ be two integers. Assume that the
sequence
$\left(\frac{\lambda_{n+1}}{\lambda_{n}}\right)_{n\in \mathbf{N}}$ is bounded,
that the sequences
$\left (\frac{b_{n,n-1}}{\lambda_n}\right )_{n \in \mathbf{N}} , \left
(\frac{b_{n,n}}{\lambda_n} \right )_{n \in \mathbf{N}}$ tend to zero, and that
there exists a positive integer $k$ and
$0\leq r <k/2$ such that the sequences
$\displaystyle{\left(|b_{n,n+1}|\left (
\frac{\lambda_{n+1}^k}{\lambda_n^k}-1 \right )\right)_{n\in \mathbf{N}}}$,
$\left(\frac{b_{n,n}}{|\lambda_n|^r}\right)_{n\in \mathbf{N}}$ and
$\left(\frac{b_{n,n-1}}{|\lambda_n|^r}\right)_{n\in \mathbf{N}}$ are
bounded. Then for every $t\geq 0$, for every piecewise constant control
$u$,
$$ |\langle \phi_{l}, \Upsilon^{u}_t\phi_n \rangle | \leq
\frac{3^{l-n}
}{(l-n)!}  \prod_{j=l+1}^{2l-n}L(j)\left ( \int_0^t |u(\tau)|\mathrm{d}\tau   \right
)^{l-n},
$$
where for $j\in{\mathbf N}$, $L(j) = \sup_{p,q \leq j} |b_{p,q}|$.
\end{lemma}}
{
\begin{IEEEproof}
Let $K > 0$. We prove the result for $u$ piecewise constant of $L^1$-norm
smaller than $K$. By Propositions \ref{PROP_Stability_Domain} and \ref{PRO_diagonal_s_weakly_coupled}, $(A,B)$ is $k$-weakly-coupled. For every $\varepsilon >0$
by Theorem~\ref{PRO_Good_Galerkin_approximation} there exists
$N=N(K,\varepsilon) > l$ such that
$
\| \Upsilon^{u}_{t}(\phi_{n}) - X^{u}_{(N)}(t,0) \phi_{n}\| < \varepsilon
$
for every $t \geq 0$.

Consider the solution $\psi: t \mapsto X^u_{(N)}(t,0) \phi_{n}$ of
\eqref{eq:sigma} with initial state $\phi_n$.
Then $\psi(t) = e^{tA^{(N)}} \phi_n + \int_{0}^t e^{(t-s)A^{(N)}} u(s)B^{(N)}
\psi(s) \mathrm{d}s$. After $l-n$ interations we get
\begin{eqnarray*}
\lefteqn{\psi(t) = e^{tA^{(N)}} \big( \phi_n +} \\
&&  + \sum_{j=1}^{l-n-1}
\int_{0\leq s_j \leq \cdots \leq s_1 \leq t}\!\!\!\!\!\!\!\!\!\!\!\!\!\!\!\!\!\!\!\!\!\!\!\!\!\!
e^{(t-s_{1}) A^{(N)}} B^{(N)}
\cdots  e^{(s_{j-1} -s_{j})A^{(N)}} B^{(N)}\times\\
&&
\times e^{s_{j}A^{(N)}} \phi_n
\prod_{m=1}^j u(s_m)  \mathrm{d}s_1 \ldots \mathrm{d}s_j +  \\
&&  +
\int_{0\leq s_{l-n} \leq \cdots \leq s_1 \leq t}\!\!\!\!\!\!\!\!\!\!\!\!\!\!\!\!\!\!\!\!\!\!\!\!\!\!\!\!\!
e^{(t-s_{1}) A^{(\!N\!)}} B^{(N)} e^{(s_{1}-s_{2}) A^{(N)}} B^{(N)} \times\cdots\\
&&\times e^{(s_{l-n-1} -s_{l-n})A^{(\!N\!)}}\!\!B^{(\!N\!)}
\psi(s_{l-n})
\!\!\prod_{m=1}^{l-n} \!\!u(s_m)
 \mathrm{d}s_1 \ldots \mathrm{d}s_{l-n}  \!\big).
\end{eqnarray*}
For the
tri-diagonal structure of the system we have
$$
\langle \phi_{l} ,
e^{(t-s_{1}) A^{(N)}} B^{(N)}
\cdots  e^{(s_{j-1} -s_{j})A^{(N)}} B^{(N)}
e^{s_{j}A^{(N)}} \phi_n \rangle = 0
$$
for every  $0\leq s_j \leq
\cdots \leq s_1 \leq t$
 and $j \leq l-n-1$.
Then
\begin{eqnarray*}
\lefteqn{\langle \phi_{l} , \psi(t) \rangle =e^{tA^{(N)}} \times} \\ 
&&  
\int_{0\leq s_{l-n} \leq \cdots \leq s_1 \leq t}\!\!\!\!\!\!\!\!\!\!\!\!\!\!\!\!\!\!\!\!\!\!\!\!\!\!\!\!\!\!\!
 \langle \phi_{l},
e^{(t-s_{1}) A^{(N)}} B^{(N)} e^{(s_{1}-s_{2}) A^{(N)}} B^{(N)} \times\\
  &&\cdots 
\times e^{(s_{l-n-1} -s_{l-n})A^{(N)}} B^{(N)}
\psi(s_{l-n})
\rangle
\prod_{m=1}^{l-n} u(s_m) \\
&&
 \mathrm{d}s_1 \ldots \mathrm{d}s_{l-n}
\end{eqnarray*}
Now,
\begin{eqnarray}\label{eq:3k}
\lefteqn{\sup_{s_1,\ldots, s_{l-n} \in [0, t]} \| B^{(N)} e^{(s_{l-n}-s_{l-n-1})
A^{(N)}} B^{(N)}  \times\cdots} \nonumber\\
&&\times e^{(s_{2}-s_{1}) A^{(\!N\!)}}\!\!
 B^{(N)}e^{(s_{1}-t) A^{(N)}} \phi_{l}\|
\leq 3^{l-n} \!\!\!\! \prod_{j=l+1}^{2l-n}\!\!\! L(j).
\end{eqnarray}
Then
\begin{eqnarray*}
\lefteqn{|\langle \phi_{l}, \psi(t) \rangle| }\\
&\leq & 3^{l-n}\prod_{j=l+1}^{2l-n}\!\!\!
L(j)\!
\int_{0\leq s_1 \leq \cdots \leq s_{l-n} \leq t}
\!\!\!\!\!\!\!\!\!\!\!\!\!\!\!\!\!\!\!\!\!\!\!\!\!\!
|u(s_1)| \ldots |u(s_{l-n})| 
\mathrm{d}s_1 \ldots \mathrm{d}s_{l-n}\\ 
&=&
  3^{l-n} \frac{\left( \int_0^t |u(s)|\mathrm{d}s\right)^{l-n}}{(l-n)!} \prod_{j=l+1}^{2l-n} L(j) ,
\end{eqnarray*}
hence
$\displaystyle{
|\langle \phi_{l}, \Upsilon^u_t(\phi_n) \rangle| \leq 3^{l-n}  \frac{K^{l-n}}{(l-n)!}
\prod_{j=l+1}^{2l-n} L(j)+ \varepsilon}$,
and the result follows as $\varepsilon$ tends to zero.
\end{IEEEproof}
}

From a physical point of view, Lemma \ref{PRO_fake_Taylor}
provides an
estimation of the probability of energy transitions (in the spirit, for
instance, of \cite[Section X.12, Example 1]{reed-simon-2}).

\begin{remark}\label{rk:2k}
In the case in which the diagonal of $B$ is zero then
equation~\eqref{eq:3k} reads
\begin{eqnarray*}
\lefteqn{\sup_{s_1,\ldots, s_{l-n} \in [0, t]} \| B^{(N)} e^{(s_{l-n}-s_{l-n-1}) A^{(N)}} B^{(N)}  \cdots} \nonumber\\
&&e^{(s_{2}-s_{1}) A^{(N)}}
 B^{(N)}e^{(s_{1}-t) A^{(N)}} \phi_{l}\|
\leq 2^{l-n}\prod_{j=l+1}^{2l-n}\!\!\! L(j).
\end{eqnarray*}
 This gives the better estimate 
$ |\langle \phi_{l}, \Upsilon^{u}_t\phi_1 \rangle | \leq
2^{l-1}
  \prod_{j=l+1}^{2l}L(j)\left ( \int_0^t |u(\tau)| {\mathrm{d}}\tau   \right
)^{l-1}/(l-1)!.
$
\end{remark}

\begin{proposition}\label{PRO_GGA_explicit_triangulaire}
 Let $(A,B)$ be a tri-diagonal system and $l$ be an integer. Assume that the
sequence
$\left(\frac{\lambda_{n+1}}{\lambda_{n}}\right)_{n\in \mathbf{N}}$ is bounded,
that the sequences
$\left (\frac{b_{n,n-1}}{\lambda_n}\right )_{n \in \mathbf{N}} , \left
(\frac{b_{n,n}}{\lambda_n} \right )_{n \in \mathbf{N}}$ tend to zero, and that
there exists a positive integer $k$ and $0\leq   r <k/ 2$ such that the sequences      
$\displaystyle{\left(|b_{n,n+1}|\left (
\frac{\lambda_{n+1}^k}{\lambda_n^k}-1 \right )\right)_{n\in \mathbf{N}}}$,
$\left(\frac{b_{n,n}}{|\lambda_n|^r}\right)_{n\in \mathbf{N}}$ and
$\left(\frac{b_{n,n-1}}{|\lambda_n|^r}\right)_{n\in \mathbf{N}}$ are
bounded. Then for every $N$ in $\mathbf{N}$, for every $t\geq 0$, $n \leq N \in \mathbf{N}$, for every piecewise constant control
$u$, 
\begin{eqnarray*}
\lefteqn{\|\pi_N \Upsilon^{u}_{t}(\phi_n) -
X^{u}_{(N)}(t,0) \phi_n\|  \leq}\\
&&
\frac{3^{N-n}
}{(N-n)!}   L(N+1)\!\!\!\! \prod_{j=N+1}^{2N-n}\!\!\!L(j)\left ( \int_0^t |u(\tau)|\mathrm{d}\tau   \right
)^{N-n+1}
\end{eqnarray*}
where for $j\in{\mathbf N}$, $L(j) = \sup_{p,q \leq j} |b_{p,q}|$.
\end{proposition}

\begin{IEEEproof}
Because of the tri-diagonal structure, (\ref{EQ_preuve_good_Galerkin}) gives
\begin{eqnarray*}
\lefteqn{\|\pi_N \Upsilon^{u}_{t}(\phi_n) -
X^{u}_{(N)}(t,0) \phi_n\| }\\
&\leq  & |b_{N,N+1}|
\left ( \int_0^t |u(\tau)|\mathrm{d}\tau \right ) \sup_{\tau \in [0,t]}|\langle \phi_{N}, \Upsilon^{u}_\tau \phi_n \rangle |.
\end{eqnarray*}
Conclusion follows from Lemma \ref{PRO_fake_Taylor}.
\end{IEEEproof}

\subsection{Example: orientation of a rotating molecule in the plane
II}\label{SEC_Rotation_Tridiagonal}
The system of Section~\ref{sec:molecule2d} provides also an example of
tri-diagonal system.
Recall that for this system, for every $j,l$ in $\mathbf{N}$, $\lambda_l=l^2$,
$\langle \phi_{j}, B \phi_{l} \rangle \neq 0$ if and only if
$|j-l|=1$ and  $\langle \phi_{j}, B \phi_{j+1} \rangle = -\mathrm{i}/2$.
We deduce a bound for the coupling constants from Proposition
\ref{PRO_diagonal_s_weakly_coupled}. For every $k$ in $\mathbf{N}$,
\begin{eqnarray*}
c_k(A,B)& \leq& \sup_{n \in \mathbf{N}} |\langle \phi_{n}, B \phi_{n+1} \rangle
| \left ( \frac{\lambda_{n+1}^k}{\lambda_n^k}-1 \right )\\
& =& \sup_{n\in
\mathbf{N}} \frac{1}{2} \left ( \left ( 1 + \frac{1}{n} \right )^{2k}-1 \right
)\\ &=& \frac{2^{2k}-1}{2}.
\end{eqnarray*}
In particular $c_1(A,B)\leq 3/2$ and, by
\eqref{EQ_estimation_N_projection}, we obtain that
$\| \pi_N \Upsilon^{u}_{t}(\phi_1) -
X^{u}_{(N)}(t,0)\pi_{N} \phi_1\| < \varepsilon$ if
$
\lambda_{N+1}=(N+1)^2 > \left(  \frac{\|u\|_{L^1}
e^{3/2\|u\|_{L^1}}}{\varepsilon} \right )^{2}$.

The tri-diagonal structure allows to obtain better estimates on $N$. From
Remark  \ref{rk:2k} and Proposition \ref{PRO_GGA_explicit_triangulaire}, we get
$$
\| X^u_{(N)}(t,0)\phi-\pi_N \Upsilon^u_t(\phi)\| \leq \frac{K^{N-1}}{(N-2)!}
$$
for any $u$ such that $\|u\|_{L^1}\leq K$ and any $\phi$ in $\mathrm{span}(\phi_1,\phi_2)$ with norm 1. 

The second estimates is significantly better than the first one. For instance,
if one has  $\|u\|_{L^1}=4$ and one desires $\varepsilon<10^{-4}$, the condition
$\varepsilon (N+1)> \|u\|_{L^1} e^{3/2\|u\|_{L^1}}$ is false for every
$N<2.7~ 10^6$   while the second condition, $\|u\|_{L^1}^{N-1}  <
{\varepsilon}(N-2)! $, is true for $N=20$.

\subsection{Example: Lyapunov design of open-loop control of the rotation of a planar molecule} \label{SEC_Lyapunoiv}
A classical method to design controls steering the system (\ref{eq:main}) from a given source to a given target is to use Lyapunov techniques (see \cite{DongPetersen}). In practice, a suitable function $V:H\rightarrow H$ is used to measure the distance between the current point $\psi$ and the target (that could be a precise wave function  or a subset of the Hilbert sphere of $H$). Under suitable regularity assumptions, the mapping $t\mapsto V(\psi(t))$ is differentiable and
$$\frac{\mathrm{d}}{\mathrm{d}t}V(\psi(t))=\mathrm{D}_{\psi(t)}V ((A+u(t) B)\psi(t))$$ is an affine function in $u(t)$.
A suitable choice of $u(t)$ depending of $\psi(t)$ ensures that the function $t\mapsto V(\psi(t))$ is decreasing. The proof that $\psi$ actually converges to the target is  non-trivial and usually relies on LaSalle invariance principles. 

When the system $(A,B)$ is weakly-coupled, the Good Galerkin Approximations may be used to obtain precise estimates on the quality of the controls obtained with Lyapunov
 techniques. For instance,  consider the system of Section~\ref{sec:molecule2d} with the source equal to $\phi=\cos(\eta)\phi_1+\sin(\eta)\phi_2$, and the target equal to 
$\{e^{i\theta} \phi_2|\theta \in \mathbf{R}\}$ where  $\phi_1$ and $\phi_2$ are the first eigenstates of the Laplacian and $\eta=10^{-3}$. On the Galerkin approximation of size 
$N=20$ of the system (\ref{eq:main}), we use the Lyapunov function $V:\psi \mapsto 1-|\langle \phi_2 , \psi \rangle|^2$ which satisfies 
$$
\frac{\mathrm{d}}{\mathrm{d}t}  V(\psi(t)) =-2u(t) \Re (\langle \phi_2,B\psi\rangle \overline{\langle \phi_2,\psi \rangle})
$$
To ensure that $V$ decreases along the trajectories of (\ref{eq:main}), we chose, for every $t$, $\tilde{u}(t):=\Re(\langle \phi_2,X^u_{(20)}(t,0)(\phi) 
\rangle$ $\langle \phi_2,BX^u_{(20)}(t,0)(\phi) \rangle)$.
We find numerically $|\langle \phi_2, X^{\tilde{u}}_{(20)}(120,0)(\cos(\eta)\phi_1+\sin(\eta)\phi_2)\rangle | >1-3.10^{-8}$ and $\int_0^{120}|u(t)|\mathrm{d}t < 4$, see~\cite{Simul} for the source of the \textbf{\texttt{Scilab}} program. From Proposition \ref{PRO_GGA_explicit_triangulaire}, we obtain  $|\langle \phi_2, \Upsilon^{\tilde{u}}_{120,0}(\cos(\eta)\phi_1+\sin(\eta)\phi_2)\rangle | >1- 10^{-4}$.

\subsection{Example: quantum harmonic oscillator}

The quantum harmonic oscillator is among the most important examples of quantum
system (see, for instance, \cite[Complement $G_V$]{cohen77}). Its controlled
version has been extensively studied (see, for instance, \cite{Rouchon,illner}).
In this example $H=L^2(\mathbf{R},\mathbf{C})$
and equation (\ref{eq:main}) reads
\begin{equation}\label{EQ_harmonic_oscillator}
\mathrm{i}\frac{\partial \psi}{\partial t}(x,t)=\frac{1}{2}(-\Delta
+x^2)\psi(x,t) +u(t) x \psi(x,t).
\end{equation}
A Hilbert basis of $H$ made of eigenvectors of $A$ is given by the sequence of
the Hermite functions $(\phi_n)_{n \in \mathbf{N}}$, associated with the
sequence $( - \mathrm{i} \lambda_n)_{n \in \mathbf{N}}$ of eigenvalues where
$\lambda_n=n-1/2$ for every $n$ in $\mathbf{N}$. In the basis $(\phi_n)_{n \in
\mathbf{N}}$, $B$ admits a tri-diagonal structure
$$
\langle \phi_j,B\phi_k\rangle = \left \{\begin{array}{cl}
- \mathrm{i} \sqrt{k-1} & \mbox{if } j=k-1\\
- \mathrm{i}\sqrt{k} & \mbox{if } j=k+1\\
0 & \mbox{otherwise}
\end{array} \right.
$$
Proposition~\ref{PROP_Stability_Domain} and
Proposition~\ref{PRO_diagonal_s_weakly_coupled}
 apply so that, for every $k$ in $\mathbf{N}$, the system $(A,B)$ is
$k$-weakly-coupled and
\begin{eqnarray*}
c_k(A,B) &\leq &\sup_n \sqrt{n} \left ( \frac{(n+1/2)^k}{(n-1/2)^k}-1 \right )
\\
&\leq & \sup_n   \sqrt{n} \left ( 1+ \frac{1}{n-\frac{1}{2}} -1 \right )\!\!
\sum_{j=0}^{k-1} \! \left ( 1+ \frac{1}{n-\frac{1}{2}} \right )^j\\
&\leq & \frac{3^{k}-1}{2} \sup_n \frac{\sqrt{n}}{n-\frac{1}{2}} \\
&\leq & 3^{k}-1.
\end{eqnarray*}

The quantum harmonic oscillator is not controllable (in any reasonable
sense) as proved in~\cite{Rouchon}.
 However, the Galerkin approximations of
(\ref{EQ_harmonic_oscillator}) of every order are exactly controllable
(see~
   \cite{PhysRevA.63.063410}), and Theorem~\ref{PRO_Good_Galerkin_approximation} ensures that any trajecory of the infinite dimensional system is a uniform limit of trajectories of its Galerkin approximations. This is not a
contradiction, since the infinite dimensional system cannot track, in general, every trajectory of its Galerkin approximations. In particular, there is no reason for which the infinite dimensional system could track a sequence of trajectories of its Galerkin approximations associated with controls with $L^1$ norm tending to infinity.
As a matter of fact, if one wants
to steer a solution of the Galerkin approximation of order $N$ of (\ref{EQ_harmonic_oscillator}) from a given state (say, the first eigenstae) to an
$\varepsilon$-neighbourhood of a given target (say, the second eigenstate), the $L^1$ norm of the control blows up as $N$ tends to infinity.

To obtain an estimate of the order $N$ of the Galerkin approximation whose
dynamics remains  $\varepsilon$ close to the one of the infinite dimensional
system when using control with $L^1$-norm $K$, one could use
Theorem~\ref{PRO_Good_Galerkin_approximation} with $k=2$, $r=1$,
$d=1$, and $\|\phi_{1}\|_{1} = 1/2$.
 The resulting bound, as given by \eqref{EQ_estimation_N_projection},
\begin{equation}\label{EQ_N_osc_harm_conservative}
N > \frac{K^{2} e^{16 K}}{4\varepsilon^{2}}  - \frac{1}{2}
\end{equation}
is however rather weak. As in the example of Section \ref{SEC_Rotation_Tridiagonal},
the tri-diagonal structure of $B$ allows better estimates.
Using Remark~\ref{rk:2k}, we find that $\|X^{(N)}_u(t,0) \phi_1 -
\pi_N \Upsilon^u_t \phi_1\|\leq \varepsilon$ provided $\|u\|_{L^1} \leq K$ and
$$
\frac{2^{N-1} \sqrt{N+2}}{(N-1)!}\sqrt{\frac{(2N)!}{(N+1)!}} K^N<\varepsilon
$$
For instance, if $K=3$ and $\varepsilon=10^{-4}$, this is true for $N= 413$, while
(\ref{EQ_N_osc_harm_conservative}) is false for $N<10^{29}$.


\section{Conclusion}

In our study we focused on the notion of weak-coupling. We  established
some interesting consequence in control theory and in numerical simulations which
applies to common physical models. We prove a result, Theorem~\ref{PRO_Good_Galerkin_approximation}, providing a uniform bound on the difference from dynamics of a finite dimensional Galerkin approximation and dynamics of the infinite dimensional system. Moreover,  an estimate on size of the Galerkin approximation has been explicitly provided
for some relevant class of systems, allowing, in particular, \emph{a priori} estimates on the error in numerical simulations on finite dimensional approximations. 
In some case, the result permits to adapt finite dimensional control techniques to study the challenging problem of the control of the bilinear Schr\"odinger equation.
For this reason we believe that the notion of weak-coupling will be a main tool in the study of controllability with relaxed controls, such as Dirac impulses, which represent, in some case, a better modelization of the physical experiences. Finally, we believe that the strong properties of convergence of the finite dimensional approximations of weakly-coupled systems will allow to address the study of a general controllabilty result for the bilinear Schr\"odinger equation with mixed spectrum.

\appendix\section{Appendices}

\subsection{Proof of Lemma~\ref{lem:Interpolation}}\label{appA}

\begin{IEEEproof}[Proof of Lemma~\ref{lem:Interpolation}]
Without loss of generality we can assume that the operators $|A|$ and $|A'|$
are positive and invertible.
 Let $(\phi_n)_{n\in{\mathbf N}}$ and $(\phi'_n)_{n\in \mathbf{N}}$
be unitary bases of $H$ made of eigenvectors
of $A$ and $A'$ respectively. Then  $\lambda_n\phi_n=|A|\phi_n$ for $n\in{\mathbf N}$ and
$D(|A|^s) = \{\psi \in
H\,:\,   \sum_{j\in \mathbf{N}} \lambda_j^{2s}
|\langle \phi_j, \psi \rangle|^{2}  < +\infty\}$. Similarly,
we can define
$\lambda'_n$  and
$D(|A'|^{s})$.

Since $D(|A|^k)\subset D(|A'|^k)$, by the closed graph theorem, we
deduce the existence of $C_k>0$ such that for every $\psi \in D(|A|^k)$
$$
\sum_n{\lambda'}_n^{2k}|\langle\psi,\phi'_n\rangle|^{2}
\leq C_k \sum_n\lambda_n^{2k}|\langle\psi, \phi_n \rangle|^{2}
$$
so that
$$
\sum_n{\lambda'}_n^{2k}\left|\sum_j\langle\psi,\phi_j\rangle\langle\phi_j,
\phi'_n\rangle\right|^{2}
\leq C_k \sum_n\lambda_n^{2k}|\langle\psi, \phi_n \rangle|^{2}.
$$
For all $\psi \in D(|A|^k)$,
let $\widetilde{\psi}$ in
$H$ such that $\psi=|A|^{-k}\widetilde{\psi}=\sum
\lambda_j^{-k}\langle\widetilde{\psi},\phi_j\rangle \phi_j$. Then, for all
$\widetilde{\psi}\in H$, we have
\begin{eqnarray}\label{eq:nabile}
\lefteqn{\sum_{n}{\lambda'}_n^{2k}
\sum_{l}\lambda_l^{-k}\overline{\langle \widetilde{\psi},
\phi_l\rangle\langle\phi_l,\phi'_n\rangle  }
\sum_{j}\lambda_j^{-k}\langle \widetilde{\psi},
\phi_j\rangle\langle\phi_j,\phi'_n\rangle} \nonumber \\
&&\leq C_k \|\widetilde{\psi}\|^{2}.~~~~~~~~~~~~~~~~~~~~~~~~~~~~~~~~~~~~~~~~~~~~~~~~~~
\end{eqnarray}
and the equality holds for $k=0$ and $C_{0}=1$.
 Consider
$\widetilde{\psi} \in
H$ and
$f_{\widetilde{\psi}}:z=s+i y\mapsto\sum_{n}{\lambda'}_n^{2(s+\mathrm{i}y)}
\langle |A|^{-s +\mathrm{i} y} \widetilde{\psi} ,\phi'_n\rangle\langle \phi'_n,
|A|^{-s -\mathrm{i} y} \widetilde{\psi} \rangle $
where, for every $z$ in $\mathbf{C}$, $|A|^{z} \widetilde{\psi} =\sum_j
\lambda_j^{z}
\langle\widetilde{\psi},\phi_j\rangle\phi_j$. Then, by \eqref{eq:nabile} for $s=0$ and $s=k$ we have
$\displaystyle{
 \left| f_{\widetilde{\psi}}(s+\mathrm{i}y)\right|\leq C_s\||A|^{-s +\mathrm{i} y}
\widetilde{\psi}\|_s \||A|^{-s -\mathrm{i} y} \widetilde{\psi} \|_s\leq C_s
\|\widetilde{\psi}\|^2}$.

If $\widetilde{\psi}$ is finite linear combination of the vectors
$\{\phi_j\}_{j\in\mathbf{N}}$ then the function $f_{\widetilde{\psi}}$
analytic on the
strip $\{z\in \mathbf{C}\,:\,0< \Re z < k\}$ and continuous on its closure
as uniform limits of a partial sum on $n$. Since it is bounded on the
boundary, by Hadamard three-lines theorem \cite[Appendix IX.4]{reed-simon-2}, it
is
bounded on the strip, and, moreover,
$ \log(\sup_{\Re z = s}|f_{\widetilde{\psi}}(z)|)$,
is a convex function of $s\in[0,k]$. So that for $s\in (0,k)$, we obtain
$\displaystyle{
\sum_n{\lambda'}_n^{2s}|\langle\psi,\phi'_n\rangle|^{2s}
\leq C_k^{\frac{s}{k}}
\sum_n|\lambda_n|^{2s}|\langle\psi, \phi_n \rangle|^{2}}
$,
and,  by density, $D(|A|^s)\subset D(|A'|^s)$. The hypothesis and the proof
being symmetric in $A$ and $A'$, we have actually the equality.
\end{IEEEproof}

\subsection{Proof of Proposition \ref{PRO_croissance_norme_A}}\label{appB}

\begin{IEEEproof}[Proof of Proposition \ref{PRO_croissance_norme_A}]
Note that for every $u\in {\mathbf R}^p$, $D(|A+\sum_l
u_lB_l|^{k/2})=D(|A|^{k/2})$, the function $|A|^{k/2}e^{t(A+\sum u_{l}
B_{l})}\psi_0$ is in $C({\mathbf R},H)$ and for every $\varepsilon>0$ the
function $|A|^{k/2}(\varepsilon( A+\sum_lu_l B_l)+1)^{-1}e^{t(A+\sum u_{l}
B_{l})}\psi_0$ is in $C^1({\mathbf R},H)$ whenever $\psi_0\in D(|A|^{k/2})$.

If $t\mapsto \psi(t)$ is the solution of \eqref{eq:main} with
 initial state $\psi_{0}$ in $D(|A|^{k/2})$, the real mapping
$f:t \mapsto
\langle |A|^k \psi(t),\psi(t) \rangle $ is absolutely continuous from ${\mathbf
R}$ to ${\mathbf R}$. We make a regularization to obtain extra
regularity, we introduce $f_\varepsilon^j:t \mapsto \langle |A|^k(\varepsilon(
A+\sum_lu_l^{j-1}B_l)+1)^{-1} \psi(t),(\varepsilon(
A+\sum_lu_l^{j-1}B_l)+1)^{-1}\psi(t) \rangle $. From the functional calculus, see~\eqref{eq:ultima} or
\cite[Theorem VIII.5]{reed-simon-1}, the sequence $f_\varepsilon^j$
is pointwise convergent to $f$ as $\varepsilon$ tends to $0$.

The function $f_\varepsilon^j$ is absolutely continuous from ${\mathbf R}$ to
${\mathbf R}$ and it is differentiable on the interval $(t_{j-1},t_j)$, for
every $t \in (t_{j-1},t_j)$,
\begin{eqnarray*}
\lefteqn{\frac{\mathrm{d}}{\mathrm{d}t} f_\varepsilon^j(t) }\\
& 
=&
 \frac{\mathrm{d}}{\mathrm{d}t}  \langle
|A|^k \Big (\varepsilon( A+\sum_lu_l^{j-1}B_l)+1 \Big )^{-1}
\psi(t),\\
&& \quad  \Big(\varepsilon( A+\sum_lu_l^{j-1}B_l)+1 \Big)^{-1}\psi(t) \rangle \\
& = & \langle |A|^k  \Big ((A+\sum_lu_l^{j-1}B_l)\varepsilon+1 \Big)^{-1} \psi(t) ,\\
&& \quad  (A+\!\!\sum_{l} u_l(t) B_l) (\varepsilon(
A+\!\!\sum_lu_l^{j-1}B_l)+1)^{-1}\psi(t)\rangle \\
&&\quad+ \langle |A|^k \Big(\varepsilon(
A+\sum_lu_l^{j-1}B_l)+1 \Big)^{-1}\\
&& \Big(A+\!\!\sum_l\! u_l(t)B_l \Big)\psi(t),  \Big (\varepsilon(
A+\!\!\sum_l\!u_l^{j-1}B_l)\!\!+\!1 \Big)^{-1}\!\!\psi(t)\rangle \\
&=& 2 \Re \langle|A|^k \Big(\varepsilon(A+\sum_lu_l^{j-1}B_l)+1 \Big)^{-1}\!\!\psi(t),\\
&& \Big(A+\sum_l u_l(t) B_l \Big) \Big(\varepsilon(
A+\sum_lu_l^{j-1}B_l)+1 \Big)^{-1}\psi(t)\rangle\\
&=& 2 \sum_l u_l(t) \Re \langle |A|^k  \Big(\varepsilon(
A+\sum_lu_l^{j-1}B_l)+1 \Big)^{-1}\psi(t),\\
&&B_l  \Big(\varepsilon(
A+\sum_lu_l^{j-1}B_l)+1 \Big)^{-1}\psi(t)\rangle,
\end{eqnarray*}
and since $(A,B_1,\ldots,B_p)$ is $k$-weakly-coupled,
\begin{eqnarray*}
\lefteqn{\left |\frac{\mathrm{d}}{\mathrm{d}t}f_\varepsilon^j(t)\right | }\\
 &\leq & 2
c_{k}(A,B_1,\ldots,B_p)\times\\
&&\quad \times \sum_l |u_l(t)| |\langle |A|^k
(\varepsilon(A+\sum_lu_l^{j-1}B_l)+1)^{-1}\psi(t),\\
&&(\varepsilon(A+\sum_lu_l^{j-1}B_l)+1)^{-1}\psi(t)
\rangle|\\
&\leq&  2 c_{k}(A,B_1,\ldots,B_p) \sum_l |u_l(t)|
f_\varepsilon^j(t).
\end{eqnarray*}
 Gronwall's lemma implies that
 $f_\varepsilon^j(t)=\langle |A|^k
(\varepsilon(A+\sum_lu_l^{j-1}B_l)+1)^{-1}\psi(t),
(\varepsilon(A+\sum_lu_l^{j-1}B_l)+1)^{-1}\psi(t) \rangle  \leq
e^{2c_{k}(A,B_1,\ldots,B_p)\sum_l \int_{t_{j-1}}^t \!|u_l|(\tau)
\mathrm{d}\tau}
f_\varepsilon^j(t_{j-1}).$
Passing to the limit $\varepsilon$ to $0$ and using \eqref{eq:ultima}, this gives
$f(t)=\langle |A|^k \psi(t),\psi(t) \rangle  \leq
e^{2c_{k}(A,B_1,\ldots,B_p)\sum_l \int_{t_{j-1}}^t \!|u_l|(\tau)
\mathrm{d}\tau}
f(t_{j-1}).$ An immediate iteration in $j$ concludes the proof.
\end{IEEEproof}

\subsection{Linear operator in Hilbert spaces}

For the reader's sake, this section recalls some basic facts of the theory of
linear operators in a Hilbert space. We refer to \cite{Kato,reed-simon-1} for
more details. 

\subsubsection{Closed operator and adjoints}
Consider a separable Hilbert space $H$
endowed with norm $\| \cdot \|$ and Hilbert product $\langle \cdot, \cdot
\rangle$. 

A linear operator is the coupled data $(A,D(A))$ where $D(A)$ is a subspace of $H$ and $A$ a
linear operator from $D(A)$ to $H$. To simplify the notation we often write $A$ instead and refer
to $D(A)$ as the domain of $A$. An operator $A'$ is an extension of $A$ if $D(A)\subset D(A')$
and $A'=A$ on $D(A)$. Below we will write $A\subset A'$. 

An operator is densely defined if its domain is dense.

An operator $A$ is closed if its graph
$
 \left\{(\psi,A\psi),\quad \psi \in D(A)\right\}
$
is a closed subspace of $H\times H$ (endowed with its natural product
topology). Notice that from the closed graph theorem, closed operator $A$ with
$D(A)=H$ are exactly bounded
operators on $H$.

An operator $A$ is closable if it has a closed extension. In this case, there
exists a smallest (in the sense of the extension) closed extension which is
called the closure and denoted $\overline{A}$. Notice that in this case the
closure of the graph of $A$ is the graph of the closure of $A$.

If $A$ is a densely defined operator, we define its adjoint $A^\ast$ by
\[
 D(A^\ast)\!=\!\left\{\phi\in H,\!\! \mbox{ s.t. } \!\exists \eta\in H,\; \!\!\forall
\psi \in D(A),\;\!
\langle \phi, A\psi \rangle \!=\!\langle \eta, \psi \rangle\right\}
\]
and for any $\phi \in D(A^\ast)$, $A^\ast\phi=\eta$, uniqueness of $\eta$
follows from the density of the domain. 

Using transformation $(\psi,\eta)\mapsto (-\eta,\psi)$ in $H\times H$,  Riesz
lemma and Closed Graph theorem we deduce that $A^\ast$ is closed and
$A$ is
closable if and only if $D(A^\ast)$ is dense, \cite[Theorem
VIII.1]{reed-simon-1}.

Notice that if $A\subset A'$ then $(A')^\ast \subset A^\ast$.

\subsubsection{Spectrum and resolvent} 
Let $A$ be a closed densely defined operator.
A complex number $\lambda$ is in the resolvent set $\rho(A)$ of $A$ if
$A-\lambda I_H$ is invertible (with bounded inverse) from $D(A)$ to $H$. The
complementary set of $\rho(A)$ is the spectrum $\sigma(A)$ of $A$.

For any $\lambda \in \rho(A)$, the operator $R_A(\lambda):=(A-\lambda
I_H)^{-1}$ is a bounded operator.  Moreover for $\lambda,\lambda' \in
\rho(A)$, $R_A(\lambda)$ commutes to $R_A(\lambda')$ and we have the following
resolvent identity
\begin{equation}\label{eq:ResolventIdentity}
R_A(\lambda)-R_A(\lambda')=(\lambda-\lambda')R_A(\lambda)R_A(\lambda').
\end{equation}
Thus for $\lambda'\neq \lambda$ in the resolvent set of $A$, we have
\[
 I_H=(\lambda'-\lambda)(I_H-(\lambda-\lambda')R_A(\lambda'))\left((\lambda'
-\lambda)^{-1}I_H -R_A(\lambda)\right)
\]
from which we deduce that the spectrum of $R_A(\lambda)$ is the closure of the
image of the spectrum of $A$ by $\lambda'\mapsto (\lambda'-\lambda)^{-1}$.

Riesz-Schauder theorem \cite[Theorem VI.15]{reed-simon-1} gives that if one of the resolvents of
$A$ is compact then the spectrum of $A$ is made of isolated eigenvalues of finite algebraic
multiplicity (the corresponding algebraic kernel is finite dimensional) possibly accumulating at
infinity. 

Notice that if one of the resolvents is compact then all of them are.

\subsubsection{Symmetric operators} A densely defined operator $A$ is symmetric if
$A\subset A^\ast$. A symmetric operator is thus always closable. A symmetric
operator is self-adjoint if $A=A^\ast$. A self-adjoint operator is thus always
closed. A symmetric operator is essentialy self-adjoint if its closure is
self-adjoint.

A densely defined operator $A$ is skew-symmetric if
$\mathrm{i}A$ is symmetric. A skew-symmetric
operator is skew-adjoint if $A=-A^\ast$, that is $\mathrm{i}A$ is self-adjoint.  A
skew-symmetric operator is essentialy skew-adjoint if its closure is
skew-adjoint.

Given a skew-adjoint operator $A$ for every $\psi$ in $D(A)$,
$$
\|(A-zI_H)\psi\|^2=\|(A+\Im zI_H)\psi\|^2+|\Re z |^2\|\psi\|^2
$$
from which we deduce that any non-purely imaginary complex number is in the
resolvent set of $A$, or the spectrum of $A$ is purely imaginary, and for every $\psi$ in $H$,
$\displaystyle{
\|R_A(z)\psi\|\leq \frac{1}{|\Re z |}\|\psi\|}$
which, from the Hille-Yosida theorem \cite[Theorem X.47a]{reed-simon-2},
provides the existence of a continuous family of unitary operators $t\in
\mathbf{R}\mapsto e^{tA}$ such that for any $\psi_0\in H$, $\psi : t\in
\mathbf{R}\mapsto e^{tA}\psi_0$ is the unique strong solution of the Cauchy
problem
$
 \partial_t \psi =A \psi \quad \psi(0)=\psi_0
$
if $\psi_0\in D(A)$ and defines a mild solution in the other cases. 

A symmetric operator $A$ is said to be positive if the associated quadratic
form $\langle A \psi,\psi\rangle$ defined on $D(A)$ is positive and it is said bounded
from below if there exists  a real $c$ such that $A-c$ is positive.

If $A$ is skew-adjoint and one of the resolvent of $A$ is compact then the
spectrum of $A$ is made of purely imaginary eigenvalues of finite multiplicity
possibly accumulating at infinity in modulus. Moreover there exists a Hilbert basis made
of eigenvectors. If $\mathrm{i}A$ is bounded from below the only accumulation point is
$+\mathrm{i} \infty$. 
 
Reciprocally if $A$ is skew-adjoint with a spectrum made of isolated
eigenvalues of finite multiplicity accumulating only at infinity then $A$ has a
compact resolvent. In this framework, the operator $A$ can be redefined the
following way (see \cite[Theorem VI.17]{reed-simon-1}). Denote by
$(\lambda_n)_{n\in\mathbf{N}}$ the spectrum of $A$ and $(\phi_n)_{n\in\mathbf{N}}$ a Hilbert basis of
eigenvectors of $A$ such that $A\phi_n=\lambda_n\phi_n$ then 
$\displaystyle{
 D(A)=\left\{\psi\in H, \quad \sum_n |\lambda_n|^2|\langle \psi,\phi_n\rangle
|^2<+\infty\right\}}$
and $ A\psi = \sum_n \lambda_n\langle \psi,\phi_n\rangle\phi_n$.

For $\psi\in D(|A|^s)$, we define
\begin{equation}\label{eq:popopo}
 |A|^s\psi = \sum_n |\lambda_n|^s\langle \psi,\phi_n\rangle\phi_n,
\end{equation}
where 
$\displaystyle{
 D(|A|^s)=\left\{\psi\in H, \quad \sum_n |\lambda_n|^{2s}|\langle \psi,\phi_n\rangle
|^2<+\infty\right\}}$.

One can also notice that for any $\psi\in H$,
$\displaystyle{
 e^{tA}\psi = \sum_n e^{t\lambda_n}\langle \psi,\phi_n\rangle\phi_n
}$
and on $D(|A|^s)$, $|A|^s$ and $e^{\mathrm{i}tA}$ commutes. Thus
\begin{equation}\label{eq:ultima}
 (\varepsilon A+1)^{-1}\psi = \sum_n (\varepsilon\lambda_n+1)^{-1}\langle
\psi,\phi_n\rangle\phi_n
\end{equation}
which, by the Dominated Convergence Theorem, tends to $\psi$  in $D(|A|^s)$, for
any $s\in\mathbf{R}$, as $\varepsilon$ goes to $0$.

\section*{Acknowledgment}

It is a pleasure for the authors to thank Ugo Boscain, Chitra Rangan, Pierre Rouchon, Mario
Sigalotti and Dominique Sugny for valuable inputs and
advices.

\bibliographystyle{IEEEtran}
\bibliography{biblio}

\begin{thebibliography}{10}
\providecommand{\url}[1]{#1}
\csname url@samestyle\endcsname
\providecommand{\newblock}{\relax}
\providecommand{\bibinfo}[2]{#2}
\providecommand{\BIBentrySTDinterwordspacing}{\spaceskip=0pt\relax}
\providecommand{\BIBentryALTinterwordstretchfactor}{4}
\providecommand{\BIBentryALTinterwordspacing}{\spaceskip=\fontdimen2\font plus
\BIBentryALTinterwordstretchfactor\fontdimen3\font minus
  \fontdimen4\font\relax}
\providecommand{\BIBforeignlanguage}[2]{{%
\expandafter\ifx\csname l@#1\endcsname\relax
\typeout{** WARNING: IEEEtran.bst: No hyphenation pattern has been}%
\typeout{** loaded for the language `#1'. Using the pattern for}%
\typeout{** the default language instead.}%
\else
\language=\csname l@#1\endcsname
\fi
#2}}
\providecommand{\BIBdecl}{\relax}
\BIBdecl

\bibitem{bms}
J.~M. Ball, J.~E. Marsden, and M.~Slemrod, ``Controllability for distributed
  bilinear systems,'' \emph{SIAM J. Control Optim.}, vol.~20, no.~4, pp.
  575--597, 1982.

\bibitem{turinici}
G.~Turinici, ``On the controllability of bilinear quantum systems,'' in
  \emph{Mathematical models and methods for ab initio Quantum Chemistry}, ser.
  Lecture Notes in Chemistry, M.~Defranceschi and C.~{Le Bris}, Eds.,
  vol.~74.\hskip 1em plus 0.5em minus 0.4em\relax Springer, 2000.

\bibitem{beauchard}
K.~Beauchard, ``Local controllability of a 1-{D} {S}chr{\"o}dinger equation,''
  \emph{J. Math. Pures Appl.}, vol.~84, no.~7, pp. 851--956, 2005.

\bibitem{camillo}
K.~Beauchard and C.~Laurent, ``Local controllability of 1{D} linear and
  nonlinear {S}chr{\"o}dinger equations with bilinear control,'' \emph{J. Math.
  Pures Appl.}, vol.~94, no.~5, pp. 520--554, 2010.

\bibitem{Nersy}
\BIBentryALTinterwordspacing
V.~Nersesyan, ``Growth of {S}obolev norms and controllability of the
  {S}chr{\"o}dinger equation,'' \emph{Comm. Math. Phys.}, vol. 290, no.~1, pp.
  371--387, 2009. [Online]. Available:
  \url{http://dx.doi.org/10.1007/s00220-009-0842-0}
\BIBentrySTDinterwordspacing

\bibitem{Schrod}
T.~Chambrion, P.~Mason, M.~Sigalotti, and U.~Boscain, ``Controllability of the
  discrete-spectrum {S}chr{\"o}dinger equation driven by an external field,''
  \emph{Ann. Inst. H. Poincar{\'e} Anal. Non Lin{\'e}aire}, vol.~26, no.~1, pp.
  329--349, 2009.

\bibitem{Schrod2}
U.~Boscain, M.~Caponigro, T.~Chambrion, and M.~Sigalotti, ``A weak spectral
  condition for the controllability of the bilinear {S}chr{\"o}dinger equation
  with application to the control of a rotating planar molecule,''
  \emph{Communications in Mathematical Physics}, vol. 311, no.~2, pp. 423--455,
  2012.

\bibitem{dalessandro-book}
D.~D'Alessandro, \emph{\BIBforeignlanguage{English}{{Introduction to quantum
  control and dynamics.}}}\hskip 1em plus 0.5em minus 0.4em\relax {Applied
  Mathematics and Nonlinear Science Series. Boca Raton, FL: Chapman,
  Hall/CRC.}, 2008.

\bibitem{book}
A.~A. Agrachev and Y.~L. Sachkov, \emph{Control theory from the geometric
  viewpoint}, ser. Encyclopaedia of Mathematical Sciences.\hskip 1em plus 0.5em
  minus 0.4em\relax Berlin: Springer-Verlag, 2004, vol.~87, control Theory and
  Optimization, II.

\bibitem{cohen77}
C.~Cohen-Tannoudji, B.~Diu, and F.~Lalo{\"e}, \emph{Quantum mechanics}, ser.
  Quantum Mechanics.\hskip 1em plus 0.5em minus 0.4em\relax Wiley, 1977.

\bibitem{ACCFEPS}
N.~Boussa{\"i}d, M.~Caponigro, and T.~Chambrion, ``Periodic control laws for
  bilinear quantum systems with discrete spectrum,'' in \emph{Proceedings of
  the American Control Conference}, 2012, pp. 5619--5824.

\bibitem{Time}
------, ``Small time reachable set of bilinear quantum systems,'' in
  \emph{Proceedings of the 51st Conference on Decision and Control (CDC)},
  2012, pp. 1083--1087.

\bibitem{ACCgates}
------, ``Implementation of logical gates on infinite dimensional quantum
  oscillators,'' in \emph{Proceedings of the American Control Conference},
  2012, pp. 5825--5830.

\bibitem{eberly-law}
C.~K. Law and J.~H. Eberly, ``Arbitrary control of a quantum electro-magnetic
  field,'' \emph{Phys. Rev. Lett.}, vol.~76, no.~7, pp. 1055--1058, 1996.

\bibitem{Boscain_Adami}
R.~Adami and U.~Boscain, ``Controllability of the schroedinger equation via
  intersection of eigenvalues,'' in \emph{Proceedings of the 44th IEEE
  Conference on Decision and Control, December 12-15}, 2005, pp. 1080--1085.

\bibitem{puel}
S.~Ervedoza and J.-P. Puel, ``Approximate controllability for a system of
  {S}chr{\"o}dinger equations modeling a single trapped ion,'' \emph{Ann. Inst.
  H. Poincar{\'e} Anal. Non Lin{\'e}aire}, vol.~26, no.~6, pp. 2111--2136,
  2009.

\bibitem{rangan}
A.~M. Bloch, R.~W. Brockett, and C.~Rangan, ``Finite controllability of
  infinite-dimensional quantum systems,'' \emph{IEEE Trans. Automat. Control},
  vol.~55, no.~8, pp. 1797--1805, 2010.

\bibitem{reed-simon-4}
M.~Reed and B.~Simon, \emph{Methods of modern mathematical physics. {IV}.
  {A}nalysis of operators}.\hskip 1em plus 0.5em minus 0.4em\relax New York:
  Academic Press [Harcourt Brace Jovanovich Publishers], 1978.

\bibitem{reed-simon-2}
------, \emph{Methods of modern mathematical physics. {II}. {F}ourier analysis,
  self-adjointness}.\hskip 1em plus 0.5em minus 0.4em\relax New York: Academic
  Press [Harcourt Brace Jovanovich Publishers], 1975.

\bibitem{buser}
P.~Buser, \emph{Geometry and spectra of compact {R}iemann surfaces}, ser.
  Modern Birkh{\"a}user Classics.\hskip 1em plus 0.5em minus 0.4em\relax
  Boston, MA: Birkh{\"a}user Boston Inc., 2010, reprint of the 1992 edition.

\bibitem{gilbarg}
D.~Gilbarg and N.~S. Trudinger, \emph{Elliptic partial differential equations
  of second order}, ser. Classics in Mathematics.\hskip 1em plus 0.5em minus
  0.4em\relax Berlin: Springer-Verlag, 2001, reprint of the 1998 edition.

\bibitem{davies}
E.~B. Davies, \emph{Spectral theory and differential operators}, ser. Cambridge
  Studies in Advanced Mathematics.\hskip 1em plus 0.5em minus 0.4em\relax
  Cambridge: Cambridge University Press, 1995, vol.~42.

\bibitem{DongPetersen}
D.~Dong and I.~Petersen, ``Quantum control theory and applications: A survey,''
  \emph{IET Control Theory \& Applications}, no.~12, pp. 2651--2671, 2010.

\bibitem{Simul}
\BIBentryALTinterwordspacing
Scilab file complement to ``{W}eakly-coupled systems in quantum control''.
  [Online]. Available: \url{http://hal.archives-ouvertes.fr/hal-00620733/en}
\BIBentrySTDinterwordspacing

\bibitem{Rouchon}
M.~Mirrahimi and P.~Rouchon, ``Controllability of quantum harmonic
  oscillators,'' \emph{IEEE Trans. Automat. Control}, vol.~49, no.~5, pp.
  745--747, 2004.

\bibitem{illner}
R.~Illner, H.~Lange, and H.~Teismann, ``Limitations on the control of
  {S}chr{\"o}dinger equations,'' \emph{ESAIM Control Optim. Calc. Var.},
  vol.~12, no.~4, pp. 615--635 (electronic), 2006.

\bibitem{PhysRevA.63.063410}
\BIBentryALTinterwordspacing
S.~G. Schirmer, H.~Fu, and A.~I. Solomon, ``Complete controllability of quantum
  systems,'' \emph{Phys. Rev. A}, vol.~63, p. 063410, May 2001. [Online].
  Available: \url{http://link.aps.org/doi/10.1103/PhysRevA.63.063410}
\BIBentrySTDinterwordspacing

\bibitem{reed-simon-1}
M.~Reed and B.~Simon, \emph{Methods of modern mathematical physics. {I}.
  {F}unctional analysis}.\hskip 1em plus 0.5em minus 0.4em\relax New York:
  Academic Press, 1972.

\bibitem{Kato}
T.~Kato, \emph{Perturbation theory for linear operators}, ser. Classics in
  Mathematics.\hskip 1em plus 0.5em minus 0.4em\relax Berlin: Springer-Verlag,
  1995, reprint of the 1980 edition.

\end{thebibliography}

%

\begin{biography}[{\includegraphics[width=1in,height=1.25in,clip, keepaspectratio]{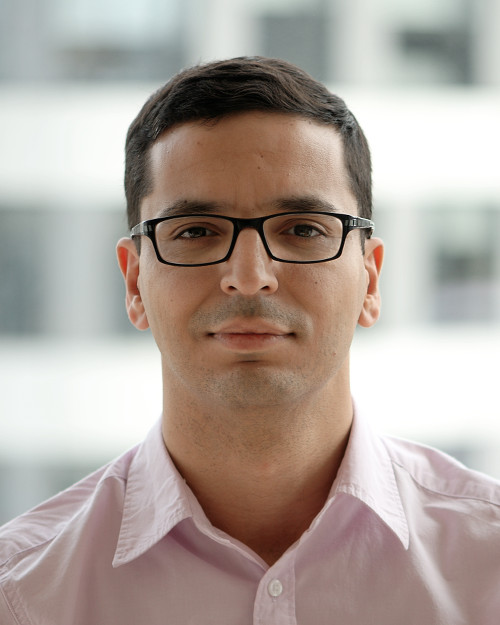}}]{Nabile Boussa\"{i}d}
  was born in France in 1978. He received his Ph.D. Degree in
applied mathematics from the University of Paris-Dauphine,  France, in 2006.

From 2002 to 2005, he was a Lecturer in mathematics at the University of Cergy-Pontoise, 
France. From 2005 to 2006, he was a Lecturer in mathematics at the University of 
Paris-Dauphine, France. From 2006 to 2007, he was a Research Associate in mathematics at Heriot-Watt University, Edinburgh, United-Kingdom.
Since 2007, he has been Assistant Professor at the University of Franche-Comt{\'e}, 
Besan\c{c}on, France. His current research interests include the application of spectral 
theory to the nonlinear partial differential equations.
\end{biography}

\begin{biography}[{\includegraphics[width=1in,height=1.25in,clip, keepaspectratio]{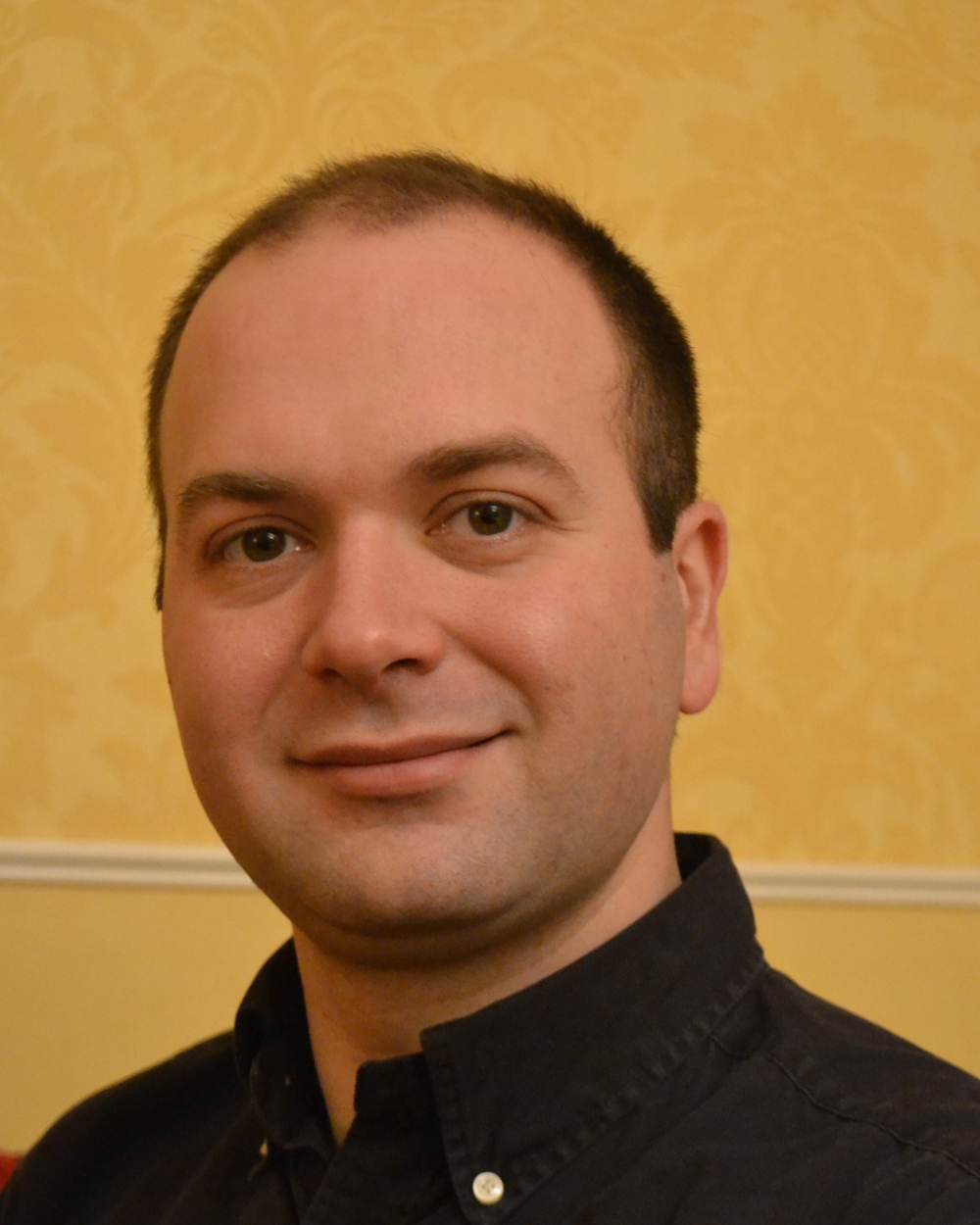}}]{Marco Caponigro}
was born in Italy in 1983. He received his Ph.D. Degree in Applied Mathematics in 2009 at SISSA, Trieste.

From 2010 to 2011 he has been postdoctoral research associate  at INRIA, Nancy Grand-Est, France. 
From 2011 to 2012 he has been postdoctoral research associate  at Center for Computation and Integrative Biology, Rutgers University, NJ, USA. 
Since 2012 he is Assistant Professor at Conservatoire National des Arts et M\'etiers, Paris, France.
\end{biography}

\begin{biography}[{\includegraphics[width=1in,height=1.25in,clip, keepaspectratio]{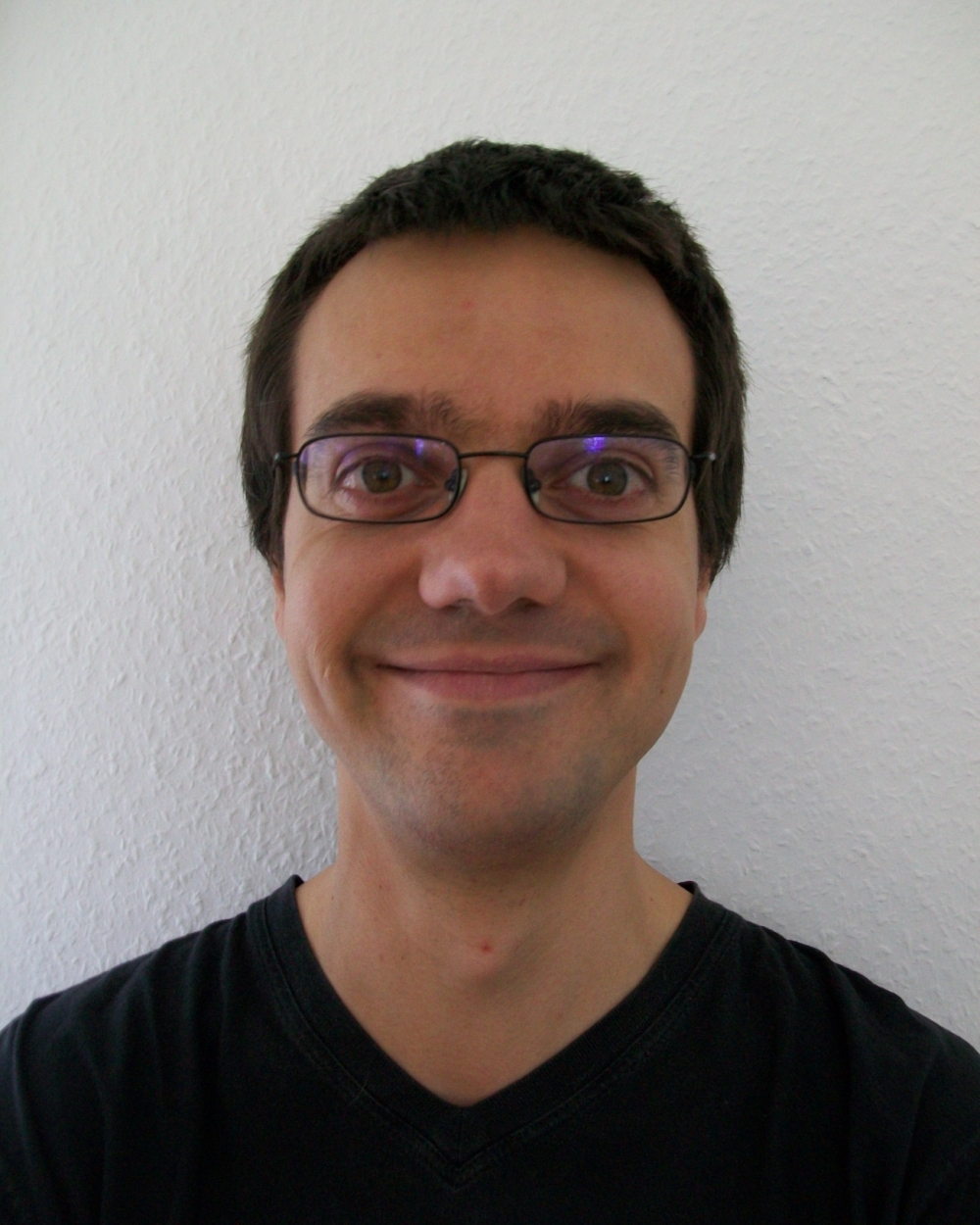}}]{Thomas Chambrion}
 was born in France in 1977. He received the  Ph.D. Degree in applied mathematics from the University of Burgundy, Dijon, France, in 2004. 

From 2001 to 2004, he was a Lecturer in mathematics at the University of Burgundy, Dijon, France. 
Since 2005, he has been teaching applied mathematics at ESSTIN, Nancy, France. His current research interests include the application of 
geometric control theory to the control of infinite dimensional systems. 
\end{biography}








\end{document}